\documentclass[11pt,leqno]{amsart}

\usepackage[utf8x]{inputenc}
\usepackage[T1]{fontenc}
\usepackage[english]{babel}

\usepackage{amsmath,amssymb,amsthm}
\usepackage{amsrefs}
\usepackage[dvipsnames]{xcolor}
\usepackage[colorlinks=true, linkcolor=MidnightBlue, citecolor=OliveGreen]{hyperref}

\usepackage{tikz}
\usetikzlibrary{automata}

\theoremstyle{definition}
\newtheorem{definition}{{\bf Definition}}[section]
\newtheorem{eg}[definition]{{\bf Example}}
\newtheorem{rmk}[definition]{{\bf Remark}}
\theoremstyle{theorem}
\newtheorem{lemma}[definition]{{\bf Lemma}}
\newtheorem{theorem}{{\bf Theorem}}
\newtheorem*{theorem*}{Theorem}
\newtheorem{cor}[definition]{{\bf Corollary}}
\newtheorem{qst}[definition]{{\bf Question}}
\newtheorem{prop}[definition]{{\bf Proposition}}

\DeclareMathOperator{\Z}{\mathbb{Z}}
\DeclareMathOperator{\N}{\mathbb{N}}
\DeclareMathOperator{\id}{id}
\newcommand*{\defeq}{\mathrel{\vcenter{\baselineskip0.5ex \lineskiplimit0pt
                     \hbox{\scriptsize.}\hbox{\scriptsize.}}}%
                     =}
					 
\DeclareMathOperator{\Sym}{\mathrm{Sym}}
\DeclareMathOperator{\Alt}{\mathrm{Alt}}
\DeclareMathOperator{\ord}{ord}
\DeclareMathOperator{\orb}{orb}
\DeclareMathOperator{\Aut}{Aut}
\DeclareMathOperator{\st}{st}
\DeclareMathOperator{\St}{St}
\newcommand{\la}{\langle}
\newcommand{\ra}{\rangle}
\newcommand{\cref}[3][]{\hyperref[#3]{#2~\ref*{#3}#1}}

\DeclareMathOperator{\syl}{syl}
\newcommand{\GGS}{\textsc{GGS} }
\newcommand{\csp}{\textsc{CS} }

\makeatletter
\@namedef{subjclassname@2020}{\textup{2020} Mathematics Subject Classification}
\makeatother

\title{Two periodicity conditions for spinal groups}

\author[J. M. Petschick]{Jan Moritz Petschick}
\address{Jan Moritz Petschick: Mathematisches
  Institut, Heinrich-Heine-Universit\"at, 40225 D\"usseldorf, Germany}
\email{jan.petschick@hhu.de}

\thanks{The research was funded by the Deutsche Forschungsgemeinschaft (DFG, German Research Foundation) — 380258175}
\date{\today}

\begin{document}

\begin{abstract}
	
	A constant spinal group is a subgroup of the automorphism group of a regular rooted tree, generated by a group of rooted automorphisms $A$ and a group of directed automorphisms $B$ whose action on a subtree is equal to the global action. We provide two conditions in terms of certain dynamical systems determined by $A$ and $B$ for constant spinal groups to be periodic, generalising previous results on Grigorchuk--Gupta--Sidki groups and other related constructions. This allows us to provide various new examples of finitely generated infinite periodic groups.
\end{abstract}

\maketitle

\section{Introduction} 
\label{sec:introduction}

Spinal groups are subgroups of the automorphism group of a regular rooted tree, that are generated by a group of rooted automorphisms $A$ and a group of directed automorphisms $B$.
Their definition is motivated by early constructions of Grigorchuk~\cite{Gri80} and Gupta and Sidki~\cite{GS83}, and they provide examples of groups of intermediate word growth, of just-infinite groups, and of finitely generated infinite periodic groups. All spinal groups are infinite, and it is easy to recognise when they are finitely generated, but it is a complicated task to obtain conditions in terms of the defining data, i.e.\ the rooted and directed groups, that ensure that a spinal group is periodic. In this paper, we establish two such conditions.

The subclasses of spinal groups for which periodicity conditions are known fall in two categories. They are either generalisations of the (first) Grigorchuk group, or of the Gupta--Sidki $p$-groups. We shall concentrate on the latter case, and consider what we will call constant spinal groups (in short \csp groups), i.e.\ spinal groups with all directed automorphisms $b$ refering directly to themselves, hence allowing a wreath decomposition
\[
	b = (b, a_1, \dots, a_{m-1}),
\]
where the elements $a_i$ are rooted. The class of \csp groups includes all groups defined by special decoration functions described in~\cite{GS84}, all (multi-)\GGS groups, and the generalised \GGS groups of Bartholdi~\cite{Bar00}. Within these subclasses, some conditions for periodicity are known, coming in two strands: a precise criterion for generalised \GGS groups with abelian rooted group given by Bartholdi~\cite{Bar00}, adapted from Vovkivski's criterion for (some) \GGS groups~\cite{Vov00}, which is related to a variant for multi-GGS groups studied by Alexoudas, Klopsch and Thillaisundaram~\cite{AKT16}. On the other strand we have a sufficient condition for groups defined by special decoration functions formulated by Gupta and Sidki~\cite{GS84}. We further develop both strands and provide two conditions for periodicity, that can be easily applied also to groups with non-cyclic directed groups and non-regular action of the rooted group. Furthermore, our results extend to non-locally finite trees and infinite rooted groups.

This is a previously unobserved phenomenon and allows us to give new examples of finitely generated infinite periodic groups. In particular, we show that the group of automorphisms of a regular rooted tree acting locally like a prescribed transitive permutation group $H$ allows infinite finitely generated periodic groups that replicate the action of $H$ on all vertices for many (but not all) permutation groups $H$.

We state our conditions in the form of certain dynamical systems on the power set of the rooted group. The system $\Lambda_b$ models the stabilised sections of automorphisms $ba$, for $a \in A$, which are itself a new tool tailor-made to observe periodic elements, see \cref{Defintion}{def:stab_sect}. The system $\Sigma_S$ models the possible local actions of automorphisms that do not reduce a certain length under taking stabilised sections. See \cref{Subsection}{sub:multi_systems} for a precise definition of the dynamical systems $\Lambda_b$ and $\Sigma_S$. We hope that these tools prove to be useful for further considerations of periodic groups within $\Aut(X^*)$. Let us now state both conditions as theorems.


\begin{theorem}\label{thm:A}
	Let $G = \la A \cup B \ra$ be a stable and strongly orbitwise-abelian \csp group with periodic rooted group $A$ such that either (i) $A$ is of finite exponent, or (ii) the directed group $B$ has finite support. Assume that $B$ is abelian and periodic. If the dynamical system $\Lambda_b$ is eventually trivial for all $b \in B$, then $G$ is periodic.
\end{theorem}

\begin{theorem}\label{thm:B}
	Let $G = \la A \cup B \ra$ be a \csp group with periodic rooted group $A$ such that either (i) $A$ is of finite exponent, or (ii) the directed group $B$ has finite support. Let $S$ be a generating set for $B$. If the dynamical system $\Sigma_S$ is eventually trivial, then $G$ is periodic.
\end{theorem}

For a \csp group to be strongly orbitwise-abelian the sections of directed elements along certain orbits of the rooted group must commute, generalising the case when the rooted group is abelian. Stability is a related homogeneity condition also trivally fulfilled by \csp groups with abelian rooted group, for precise definitions see \cref{Definitions}{def:stable} and \ref{def:orbitwise_ab}. However, \cref{Theorem}{thm:A} applies to many more groups, cf.\ \cref{Example}{eg:non-abelian strongly orbitwise-abelian stable}.

In general, the dynamical systems $\Lambda_b$ are more likely to be eventually trivial than the systems $\Sigma_S$, which justifies the additional conditions of \cref{Theorem}{thm:A}.

All previously mentioned conditions for periodicity of subclasses of \csp groups can be derived from the two theorems above; see \cref{Subsection}{sub:discussion_of_cref_theorem_thm_abelian_case} and \cref{Remark}{rmk:regular case}. Although \cref{Theorem}{thm:A} and \cref{Theorem}{thm:B} give only sufficient conditions for periodicity, in the case of an abelian rooted group the sufficient and necessary criterion of Bartholdi may be extended to a broader class of groups; cf.\ \cref{Corollary}{cor:abelian_case}.

Recently, Rajeev and the author~\cite{PR21} described a method of constructing, given a group $G \leq \Aut(X^*)$ and a positive integer $s$, a new group $\operatorname{Bas}_s(G) \leq \Aut(X^*)$ called the $s$\textsuperscript{th} Basilica group of $G$, based on the Basilica group defined by Grigorchuk and \.{Z}uk~\cite{GZ02}. It is a remarkable feature that the Basilica groups of \csp groups are \csp groups. We prove that Basilica groups of \csp groups satisfying the conditions of \cref{Theorem}{thm:B} again satisfy \cref{Theorem}{thm:B}, which is relevant to the question if Basilica groups of periodic groups are periodic in general. Combined with our good understanding of periodic \csp groups in the abelian case, this provides a wealth of examples of periodic \csp groups with non-regular rooted action. 

Both \cref{Theorem}{thm:A} and \ref{thm:B} only provide sufficient conditions, and, in addition to a family of groups that are periodic in accordance with \cref{Theorem}{thm:B} (\cref{Example}{eg:B}), we give an example of a periodic \csp group that satisfies neither the condition of \cref{Theorem}{thm:A} nor those of \cref{Theorem}{thm:B}.

\subsection*{Acknowledgements}
This is part of the author's Ph.D. thesis, written under the supervision of Benjamin~Klopsch at the Heinrich-Heine-Universität~Düsseldorf. The author thanks Benjamin~Klopsch, Karthika~Rajeev and Eike~Schulte for many helpful suggestions which led to significant improvements in the exposition of this paper.


\section{Preliminaries and \csp groups} 
\label{sec:preliminaries_and_csp_groups}
	
\subsection{Automorphisms of rooted trees} 
\label{sub:automorphisms_of_rooted_trees}

Let $X$ be a non-empty set, $0 \in X$ an element and $X^*$ the free monoid over $X$, which we identify with its Cayley graph, i.e.\ a regular rooted tree of valency $|X|+1$. We call $X$ the \textit{alphabet}, its members \textit{letters} and $0$ the \textit{distinguished letter}. Observe that we do not restrict to finite sets. Denote $X\setminus\{0\}$ by $\dot X$. If $X$ is a group, we require that the distinguished letter equals the neutral element. We write \(\epsilon\) for root, i.e.\ the empty word. The \(n\)\emph{\textsuperscript{th} layer of} \(X^*\) is the set $X^n$ of vertices represented by words of length \(n\), or, equivalently, of vertices of distance $n$ to the root.

A \textit{(rooted tree) automorphism} of $X^*$ is a graph automorphism fixing the root $\epsilon$. The invariance of the root is a feature of all graph automorphisms if $X$ is finite. Since the root is fixed, every rooted tree automorphism maps the \(n\)\textsuperscript{th} layer to itself. The action of the full group of rooted tree automorphisms \(\Aut(X^*)\) is transitive on each layer. A subgroup of $\Aut(X^*)$ with this property is called \emph{spherically transitive}. The stabiliser of a word \(u\) under the action of a group \(G \leq \Aut(X^*)\) is denoted by \(\st_G(u)\), and the intersection of all stabilisers of words of length \(n\) is called the \emph{\(n\)\textsuperscript{th} layer stabiliser}, denoted \(\St_G(n)\).

Let \(a\) be an automorphism of \(X^*\) and let \(u, v\) be words. We denote \textit{concatenation}, i.e.\ the monoid multiplication in $X^*$, by $\star$. Since layers are invariant under \(a\), the equation
\[
	(u \star v).a = u.a \star v.(a|_u)
\]
defines a unique automorphism \(a|_u\) of \(X^*\) called the \emph{section of \(a\) at \(u\)}.

Consequently, any automorphism $a$ can be decomposed into the sections prescribing the action at the subtrees of the first layer, and $a|^\epsilon$, the action of $a$ on the first layer $X^1 = X$. We adopt the convention that an $X$-indexed family of automorphisms $(a_x)_{x \in X}$ is identified with the automorphism having section $a_x$ at $x$ and $(a_x)_{x \in X}|^\epsilon = \id$. Hence for any $a$ we write
\[
	a = (a|_x)_{x \in X}\; a|^\epsilon.
\]
Conversely, any family $(a|_x)_{x \in X}$ of automorphisms together a permutation $a|^\epsilon \in \Sym(X)$ defines a unique automorphism of $X^*$.
For all words \(u, v \in X^*\) and automorphisms \(a, b \in \Aut(X^*)\) we have
\[\begin{array}{ccc}
	(a|_u)|_v = a|_{u \star v}, &
	(ab)|_u = a|_ub|_{u.a}, &
	a^{-1}|_u = (a|_{u.a^{-1}})^{-1}.
\end{array}\]
Let $H \leq \Sym(X)$ be a subgroup. By the computation rules above, we see that the set of automorphisms $a$ such that $a|_u|^\epsilon \in H$ for all $u \in X^*$ is a subgroup, denoted $\Gamma(H)$ and called the \textit{$H$-labeled subgroup of $\Aut(X^*)$}.

By convention, permutations of $X$ act on $X^*$ by permuting the first layer subtrees, i.e.\ we consider $\Sym(X)$ to be embedded into $\Aut(X^*)$ by $\rho \mapsto (\id)_{x \in X}\; \rho$. Automorphisms of this kind are called \emph{rooted}. Consequently, we think of $\Aut(X^*)$ as acting from the right, however, we conjugate from the left, writing ${}^ab = a b a^{-1}$. The purpose of this is twofold: calculation of sections of conjugates involve fewer inversions, i.e.\ for $g, h \in \Aut(X^*)$ and $u \in X^*$ we have $({}^hg)|_x = {}^{h|_x}(g|_{x.h})$, and furthermore, since on some occasions group elements or vertices will be represented with integers, we hope to better distinguish between powers and conjugation.


\subsection{Constant spinal groups} 
\label{sub:constant_spinal_groups}
	
We now define the family of constant spinal groups, which includes all (multi-)\GGS groups and the more general groups defined in~\cite{GS84} by a special decorating function. Both name and definition are derived from the family of spinal groups defined by Bartholdi and \v{S}uni\'{k}~\cites{BS01,BGS03}, which we shall not recall in detail, but we remark that the defining data of a general spinal group contains a certain sequence of homomorphisms, and the constant spinal groups are those where the defining sequence is constant.

\begin{definition}\label{def:GGS-group}
	Let $A \leq \Sym(X)$ be a transitive permutation group, embedded into $\Aut(X^*)$ as rooted automorphisms, and let $B \leq \St(1)$ be a subgroup such that
	\[
		\la B|_x \mid x \in \dot X \ra = A \quad\text{ and }\quad b|_0 = b\; \text{ for all } b \in B.
	\]
	We call the group $G = \la A, B \ra$ a \emph{constant spinal group}, or in short a \emph{\csp group}.
	The groups $A$ and $B$ are referred to as the \emph{rooted group} and the \emph{directed group of $G$}, respectively. We say that $B$ has \emph{finite support} if for all $b \in B$ and all almost all $x \in X$ we have $b|_x = \id$.
\end{definition}

A \GGS group is most commonly defined by a so-called \textit{defining vector}: Let $m \in \N_{\geq 2}$, let $X = \{0, \dots, m-1\}$, let $a \in \Sym(X)$ be an $m$-cycle, and $\underline{e} = (e_1, \dots, e_{m-1})$ a vector in the free module $(\Z/m\Z)^{m-1}$ such that $\la e_i \mid i \in \dot X\ra = \Z/m\Z$. Define an automorphism $b = (b, a^{e_1}, \dots, a^{e_{m-1}})$. Then $\la a, b \ra$ is the \GGS group with defining vector $\underline{e}$.
In the terminology of \csp groups, this is the \csp group with rooted group $\la a \ra$ and directed group $\la b \ra$.

For comparision with previous results on periodicity of spinal groups, we give a small dictionary that describes certain well-studied subfamilies of spinal groups in the terminology of \csp groups.

\begin{itemize}
	\item \GGS groups are \csp groups whose rooted group and directed group are cyclic.
	\item The multi-\GGS groups of~\cite{AKT16} are \csp groups with cyclic rooted group.
	\item The generalised \GGS groups of~\cite{Bar00} are \csp groups with cyclic directed group.
	\item The groups of~\cite{GS84} defined by a special decorating function are \csp groups with rooted group acting regularly and cyclic directed group.
	\item Neither the so-called $\mathcal G$ groups of~\cite{BS01} nor the \textsc{EGS} groups of~\cite{Per07} (both are classes of spinal groups where periodicity criteria are known) are \csp groups. This is clear for their usual embedding into $\Aut(X^*)$, and follows from~\cite{PT21} for their isomorphism classes.
\end{itemize}

It is evident that all \csp groups are \emph{self-similar}, i.e.\ that for all $x \in X$ and $g \in G$ the section $g|_x$ is contained in $G$. Furthermore they are \emph{fractal}, i.e.\ $\St_G(1)|_x = G$ for all $x \in X$: observe that ${}^ab|_x = b|_{x.a}$ with $a \in A, b \in B$ produces a generating set for $A$, as well as for $B$, using the transitivity of $A$. This implies that $G$ is spherically transitive. Clearly $G$ is contained in the $A$-labelled subgroup of $\Aut(X^*)$. Since $\St_G(1)$ is a proper subgroup with a surjection to the full group, $G$ is infinite.

Every \csp group $G$ is a homomorphic image of the free product $A \ast B$. Thus every element $g \in G$ can be written in the form
\begin{align*}
	g &= ({}^{a_0}b_0) ({}^{a_1}b_1) \dots ({}^{a_{n-1}}b_{n-1})a_n, &{a_i} \in A, b_i \in B, n \in \N.
\end{align*}
The minimum of all possible numbers $n$ in such words is called the \emph{syllable length} of $g \in G$, denoted $\syl(g)$. This is the weighted word length with respect to the generating set $A \cup B$ with elements in $A$ having weight $0$. Consequently, we refer to the elements ${}^AB$ as \emph{syllables}. Observe that all syllables stabilise the first layer, hence $g|^\epsilon = a_n$.

The common strategy for proving that a \csp group is periodic, starting from the proofs for the examples of Grigorchuk and Gupta and Sidki, is to establish a contraction of the syllable length upon taking sections. We refine the usual statements concerning this contraction by employing the new concept of stablised sections.

\begin{definition}\label{def:stab_sect}
	Let $v \in X^*$ be a vertex and $g \in \Aut(X^*)$. Let $\ell_g(v)$ be the length of the orbit of $v$ under $g$. If $\ell_g(v) < \infty$, define the \emph{stablised section of $g$ at $v$} by
	\[
		g\|_v \defeq g^{\ell_g(v)}|_v.
	\]
\end{definition}

The usefullness of stabilised sections for establishing periodicity of groups of automorphisms of $X^*$ stems from the fact that powers of $g$ stabilising the $n$\textsuperscript{th} layer can be described by certain stabilised sections. We record this and two other useful facts in a lemma.

\begin{lemma}\label{lem:power of stab sect}
	Let $g \in \Aut(X^*)$ and let $n, m \in \N$ such that $g^n \in \St(m)$. Then
	\begin{enumerate}
		\item $g^n|_v = (g\|_v)^{\frac n {\ell_g(v)}}$ for all $v \in X^m$.
	\end{enumerate}
	Let $u, v \in X^*$ such that $\ell_g(u \star v)$ is finite. We have:
	\begin{enumerate}\setcounter{enumi}{1}
		\item $\ell_g(u \star v) = \ell_g(u) \cdot \ell_{g\|_u}(v)$.
		\item $g\|_{u\star v} = g\|_u\|_v$.
	\end{enumerate}
\end{lemma}

\begin{proof}
	Since $v$ is stabilised by $g^n$, the number $n$ is a multiple of $\ell_g(v)$. Thus
	\[
		g^n|_v = g^{\ell_g(v) \cdot \frac n {\ell_g(v)}}|_v = (g\|_v)^{\frac n {\ell_g(v)}}.
	\]
	For (2), observe that $\ell_g(u\star v)$ must be a multiple of $\ell_g(u)$, and calculate
	\[
		(u \star v).g^{\ell_g(u)} = u.g^{\ell_g(u)} \star v.(g^{\ell_g(u)}|_u) = u \star v.(g\|_u).
	\]
	Using the equation again, we prove (3):
	\[
		g\|_{u \star v} = g^{\ell_g(u)\ell_{g\|_u}(v)}|_{u \star v} = (g^{\ell_g(u)}|_u)^{\ell_{g\|_u}(v)}|_v = g\|_u\|_v.\qedhere
	\] 
\end{proof}

We shall use \cref{Lemma}{lem:power of stab sect} regularly and without reference.

\begin{lemma}\label{lem:length_reduction_in_GGS}
	Let $G$ be a \csp group and let $g \in G$. We have:
	\begin{enumerate}
		\item $\syl(g|_x) \leq \frac 1 2 (\syl(g) + 1)$ for all $x \in X$.
		\item $\sum_{x \in X}\syl(g|_x) \leq \syl(g)$.
		\item $\syl(g\|_x) \leq \syl(g)$ for all $x \in X$ such that $\ell_g(x) < \infty$.
		\item If $g|^\epsilon$ has finite order, then $g\|_x \in A$ for almost all $x \in X$.
	\end{enumerate}
\end{lemma}

\begin{proof}
	(1) Consider an element $g = {}^{a_0}b_0 {}^{a_1}b_1 \dots {}^{a_{n-1}}b_{n-1}a_n \in G$, where $n = \syl(g)$. Taking the section at $x \in X$, one obtains
	\begin{equation*}\label{eq:sections}
		g|_x = b_0|_{x.a_0} b_1|_{x.a_1} \dots b_{n-1}|_{x.a_{n-1}}.\tag{$\ast$}
	\end{equation*}
	For each $i \in \{0, \dots, n-1\}$ the section $b_i|_{x.a_i}$ is an element of $G\setminus A$ if and only if $x.a_i = 0$. In this case it is a member of $B$. If two consecutive syllables both yield an element of $B$, they reduce to a single syllable. Hence at most every second syllable of $g$ produces a syllable of $g|_x$, yielding the inequality.
	
	(2) This is a consequence of the fact that every syllable produces exactly one letter from $B$ in all of its sections.
	
	(3) For any section of $g^{\ell_g(x)}$ we have
	\begin{equation*}\label{eq:orbit sections}
		g^{\ell_g(x)}|_x = g|_x g|_{x.g} \dots g|_{x.g^{\ell_g(x)-1}}.\tag{$\ast\ast$}
	\end{equation*}
	But since $\{x, x.g, \dots, x.g^{\ell_g(x)-1}\}$ is the orbit of $x$ under $g$, all sections are taken at different positions in $X$ and the assertion follows from (2).
	
	(4) If $g|^\epsilon$ has finite order, every $\la g \ra$-orbit in $X$ is finite. By (2) there are only finitely many sections outside $A$, i.e.\ there are only finitely many $x \in X$ such that $g|_x \notin A$. Thus only finitely many stabilised sections of $g$ can have positive syllable length.
\end{proof}

Since $b\|_0 = b|_0 = b$ holds for all $b \in B$, the inequality in (3) of the previous lemma cannot be strict for all non-trivial elements of a given \csp group. However, for some \csp groups it is possible to obtain a strict inequality excluding a controllable `error' set, on which our conditions for $G$ to be periodic rest.


\subsection{Generalities on periodic \csp groups} 
\label{sub:generalities_on_periodic_csp_groups}

Since we do not restrict to finite alphabets $X$, we need two lemmata to control possible problems for infinite $X$.

\begin{lemma}\label{lem:a-adic is periodic up to level}
	Let $H \leq \Sym(X)$ be periodic and $\Gamma(H) \leq \Aut(X^*)$ the $H$-labelled subgroup. Then $\Gamma(H)/\St_{\Gamma(H)}(n)$ is periodic for all $n \in \N$.
\end{lemma}

\begin{proof}
	The group $\Gamma(H)/\St_{\Gamma(H)}(n)$ is isomorphic to the $n$-fold iterated unrestricted wreath product $H \wr \dots \wr H$, which is periodic.
\end{proof}

\begin{lemma}\label{lem:finite support}
	Let $G$ be a \csp group such that the directed group $B$ has finite support, let $n \in \N$. Then for every $g \in G$ the set $\{ g|_v \mid v \in X^n \}$ is finite.
\end{lemma}

\begin{proof}
	This is an immediate consequence of the composition rule for sections.
\end{proof}

\begin{prop}\label{prop:torsion_by_length}
	Let $G$ be a \csp group with periodic rooted group $A$.
	Let either (i) $A$ be of finite exponent, or (ii) let the directed group $B$ have finite support. Let $T \subseteq G$ be a set of elements of finite order. If for every $g \in G \setminus A$ there is some $m \in \N$ such that for all $v \in X^m$ we have
	\begin{equation*}\label{eq:lr}
		\syl(g\|_v) < \syl(g)\quad \text{ or }\quad \syl(g\|_v) \in T, \tag{$\dagger$}
	\end{equation*}
	the group $G$ is periodic.
\end{prop}

\begin{proof}
	We use induction on the syllable length. Rooted elements are of finite order, thus we consider only $g$ with $\syl(g) \geq 1$.
	Let $m$ be the integer such that (\ref{eq:lr}) holds for $g$ and $X^n$. Since $G$ is contained in $\Gamma(A)$, the order of the action of $g$ on $X^m$ is a finite number $k \in \N$ by \cref{Lemma}{lem:a-adic is periodic up to level}, thus
	\[
		g^k = ((g\|_v)^{\frac k {\ell_g(v)}})_{v\in X^m}.
	\]
	By assumption, either $g\|_v$ is of shorter syllable length than $g$, or it is a member of $T$, hence of finite order. By induction $g\|_v$ is of finite order. Also, by \cref[(4)]{Lemma}{lem:length_reduction_in_GGS} at most finitely many stabilised sections have positive syllable length.

	In case $(i)$, let $k'$ be the least common multiple of the order of these sections and the exponent of $A$, in case $(ii)$ let $k'$ be the least common multiple of the orders of all (finitely many by \cref{Lemma}{lem:finite support}) non-trivial first layer sections of $g\|_v$. In both cases $(g\|_v)^{k'} = \id$ for all $v \in X^n$, hence $(g^k)^{k'} = \id$.
\end{proof}

We remark that the set $T$ has to include the directed group $B$ as a subset, since $b|_0 = b$. Therefore the existence of an error set $T$ implies that $B$ is a periodic group.


\subsection{Some dynamical systems} 
\label{sub:multi_systems}

We use \cref{Proposition}{prop:torsion_by_length} to establish our two dynamical conditions for periodicity. It is beneficial to reformulate this proposition using the language of atomic dynamical systems.

\begin{definition}
	Let $A$ be a set and let $\phi: \mathcal P (A) \to \mathcal P(A)$ be a self-map of the powerset of $A$. The dynamical system $(\mathcal P(A), \phi)$ is called \emph{atomic} if
	\[
		\phi(S) = \bigcup_{s \in S} \phi( \{ s \})
	\]
	for all $S \in \mathcal P(A)$. Clearly, an atomic dynamical system is defined by its action on singletons. It is called \emph{eventually $T$} if for all $a \in A$ there is an integer $n \in \N$ such that for all $m > n$ we have $\phi^m(\{a\}) \subset S$. If $A$ is a group, the system is \emph{eventually trivial} if it is eventually $\{ 1_A \}$.
\end{definition}

Now we may restate \cref{Proposition}{prop:torsion_by_length}. A \csp group $G$ with periodic rooted and direct group fulfilling condition (i) or (ii) as stated in the proposition is periodic if the  dynamical system $(\mathcal P (G), \|_X)$ is eventually $A \cup T$, where the self-map is atomic and hence defined by
\[
	\|_X: \mathcal P (G) \to \mathcal P (G), \quad \{g\} \mapsto \{ g\|_x \mid x \in X \}.
\]
We want to obtain, at least incase of $G$ having a finite rooted group, similar conditions on finite dynamical systems, which we shall now define. First, we need some terminology.

Let $x, y \in X$ and let $A \leq \Sym(X)$ be a group acting transitively on $X$. Write $\operatorname{mp}_A(x,y)$ for the set of elements $a \in A$ such that $x.a = y$, and fix an element $e_{x\mapsto y} \in \operatorname{mp}_A(x,y)$. Then
\[
	\operatorname{mp}_A(x,y) = e_{x\mapsto y} \cdot \st_A(y) = \st_A(x) \cdot e_{x\mapsto y}.
\]
By convention, set $e_{0 \mapsto 0} = 1_A$. The sets $\operatorname{mp}_A(0, x)$ with $x \in X$ are right cosets of the point stabiliser $\st_A(0)$, hence they form a partition of $A$, which we use to define certain subsets of the set of conjugates of a given element.

\begin{definition}\label{def:x+c}
	Let $G = \la A \cup B \ra$ be a \csp group. For $a \in A$ and $x \in X$, we define
	\[
		\mathfrak C(a, x) \defeq {}^{\operatorname{mp}_A(0, x)}a \subseteq A
		\quad\text{ and }\quad
		\mathfrak X(a, x) \defeq \bigcup_{c \in \mathfrak C(a,x)} \orb_c(0) \setminus \{ 0 \} \subseteq \dot X.
	\]
\end{definition}

For every $x$ in the orbit $\orb_a(0)$ of $0$ under $a$ we have $a \in \mathfrak C(a, x)$. If $A$ acts regularly, we identify $X$ with $A$ such that $e_{0 \to x} = x$, hence $\mathfrak C(a, x) = \{^xa\}$ and $\mathfrak X(a) = \orb_{{}^xa}(0)\setminus\{0\} = 0.(\la {}^xa \ra \setminus \{ 1_A\}$). We now define the first dynamical system.

\begin{definition}\label{def:sigma}
	Let $G = \la A \cup B \ra$ be a \csp group, let $S\subseteq B$, let $a \in A$ and let $x \in X$.
	Define a map
	\begin{equation*}\label{eg:sigma}
		\sigma_S(a, x) = \la \prod_{c' \in \la c \ra \setminus \st_A(0)} b|_{0.c'} \mid c \in \mathfrak C(a, x), b \in S \ra \cdot \la b|_y \mid y \in \mathfrak X(a, x), b \in S \ra' \subseteq A.\tag{$\mathparagraph$}
	\end{equation*}
	Notice that the order of the products generating the left side of the definition can be chosen arbitrarily, using the derived subgroup on the right side. Now the dynamical system $(\mathcal P(A), \Sigma_S)$ is the atomic system defined by
	\[
		\Sigma_S: \mathcal P(A) \to \mathcal P(A), \quad \Sigma_S(\{a\}) = \bigcup_{x \in X} \sigma_S(a, x).
	\]
\end{definition}

Our second dynamical system makes use of the stabilised section map $\|_0$.

\begin{definition}\label{def:finite_dynamics}
	Let $G = \la A \cup B \ra$ be a \csp group with periodic rooted group $A$. For every $b \in B$ we define a map $\lambda_b: A \to A$ by
	\[
		\lambda_b(a) \defeq \prod_{i = 1}^{\ell_a(0)-1} b|_{0.a^i} \quad\text{ so that }\quad (ba)\|_0 = b\lambda_{b}(a).
	\]
	Let $x \in X$. Then we define
	\[
		\lambda_b(a, x) = \lambda_b({}^{e_{0 \mapsto x}}a),
	\]
	and an atomic dynamical system $(\mathcal P(A), \Lambda_b)$ given by the map
	\[
		\Lambda_b(\{a\}) = \{ \lambda_b(a, x) \mid x \in X \}.
	\]
\end{definition}

Clearly, $\lambda_b(a, 0) = \lambda_b(a)$ for all $a \in A, b \in B$. If the group $A$ is abelian, $\lambda_b(a, x) = \lambda_b(a)$ for all $a \in A, b \in B$ and $x \in X$. Consequently, we may regard $\Lambda_b$ as a map of type $A \to A$ in the abelian case, and write $\Lambda_b(a) = \lambda_b(a)$.

\begin{definition}\label{def:stable}
	A \csp group $G$ is called \textit{stable} if for all $a \in A, x \in X$ and $b \in B$
	\[
		\lambda_b(c') = \lambda_b(c'') \text{ for all } c', c'' \in \mathfrak C(a, x).
	\]
\end{definition}

Let $G$ be a stable \csp group, and let $a, a' \in A, b \in B$. Then $\lambda_b(^{a'}a) \in \Lambda_b(a)$, since $\lambda_b(^{a'}a) = \lambda_b(^{e_{0\mapsto 0.a'}}a)$.


\subsection{(Strongly) orbitwise-abelian \csp groups} 
\label{sub:_strongly_orbitwise_abelian_groups}

Reviewing the dynamical systems on $\mathcal P(A)$ defined by $\Sigma_S$ and $\Lambda_b$ with $S \subseteq B, b \in B$, we see that (heuristically), the systems $\Lambda_b$ are more likely to be eventually trivial, since $|\Lambda_b(a)| \leq |X|$, while the case $\Sigma_S(a) = \sigma_S(a,x) = A$ does frequently (and naturally, cf.\ \cref{Subsection}{sub:limitations_and_symmetric_groups}) occur. The drawback of the dynamical systems $\Lambda_b$ is that we need some (weak) form of abelianess within $A$ to conclude the periodicity of $G$, even if all $\Lambda_b$ are eventually trivial.

\begin{definition}\label{def:orbitwise_ab}
	A \csp group is called \emph{orbitwise-abelian} if for all $a \in A$
	\[
		\la b|_{0.c} \mid c \in \la a \ra \setminus \st_A(0), b \in B \ra \leq A
	\]
	is abelian. It is called \emph{strongly orbitwise-abelian} if
	\[
		\la b|_{y} \mid y \in \mathfrak X(a, x), b \in B \ra \leq A
	\]
	is abelian for all $a \in A$ and $x \in X$.
\end{definition}

Any strongly orbitwise-abelian group is orbitwise-abelian, since $a \in \mathfrak C(a, 0) = {}^{\st_A(0)}a$, hence $0.c \in \mathfrak X(a,0)$ for all $c \in \la a \ra  \setminus \st_A(0)$. Notice that for a strongly orbitwise-abelian group, the derived group in the \cref{Definition}{def:sigma} of $\Sigma_S$ is trivial. We also gain better control over the interplay between the dynamical systems $\Lambda_b$ for various $b\in B$.

\begin{lemma}\label{lem:orbitwise_abelian_dyn}
	Let $G$ be a orbitwise-abelian \csp group with periodic rooted group. Then the map $b \mapsto \lambda_b: B \to A^A$ is a group homomorphism.
\end{lemma}

\begin{proof}
	Let $b_1, b_2 \in B$. Then for all $a \in A$,
	\begin{align*}
		\lambda_{b_1b_2}(a) = \prod_{i = 1}^{\ell_a(0)-1} (b_1b_2)|_{0.a^{i}} &= \prod_{i = 1}^{\ell_a(0)-1} b_1|_{0.a^{i}}b_2|_{0.a^{i}}\\
		&= \prod_{i = 1}^{\ell_a(0)-1} b_1|_{0.a^i}\prod_{i = 1}^{\ell_a(0)-1}b_2|_{0.a^i}
		= \lambda_{b_1}(a)\lambda_{b_2}(a).\qedhere
	\end{align*}
\end{proof}

The map $b \mapsto \lambda_b$ is not necessarly injective, even if $G$ is orbitwise-abelian. Take for example a \GGS $p$-group $G = \la a, b \ra$ with defining vector $\underline e$ on a $p$-adic rooted tree such that $\sum_{i = 1}^{p-1} \equiv_p 0$. On $p$-adic trees, these are precisely the periodic \GGS groups. We have
\[
	\lambda_b(a^j) = \prod_{i = 1}^{p-1} b|_{0.a^{ji}} = a^{\sum_{i = 1}^{p-1} e_i} = 1_A
\]
for all $j \in \{0, \dots, p-1\}$. Thus the image of $B$ under $b \mapsto \lambda_b$ is trivial.



\section{Proofs and discussion of the theorems} 
\label{sec:proofs_of_the_theorems}

\subsection{Proof of \cref{Theorem}{thm:A}} 
\label{sub:proof_of_cref_theorem_thm_abelian_case}

It is useful to consider the next lemma seperately, as it will also be used in the proof of \cref{Theorem}{thm:B}.

\begin{lemma}\label{lem:non-reduction in one step}
	Let $G = \la A \cup B \ra$ be a \csp group, let $g = {}^{a_0}b_0 \dots {}^{a_{n-1}}b_{n-1}a_n \in G$ and let $x \in X$ such that $n = \syl(g\|_x) = \syl(g)$. Then there are integers $j(i)$, for all $i \in \{0, \dots, n-1\}$, such that
	\[
		g\|_x = \prod_{j = 0}^{\ell-1} \prod_{i = 0}^{n-1} b_i|_{0.^{a_i^{-1}}(a_n^{j+j(i)})}
	\]
	where $\ell = \ell_g(x)$ and ${}^{a_i^{-1}}a_n \in \mathfrak C(a_n, x)$.
\end{lemma}

\begin{proof}
	We compute the stabilised section at $x$, see (\ref{eq:sections}) and (\ref{eq:orbit sections}).
	\begin{align*}
		g\|_x &= \left(\prod_{j = 0}^{\ell-1} \left(\left(\prod_{i = 0}^{n-1} (^{a_i}b_i) \right) a_n\right)\right)\bigg|_x = \prod_{j = 0}^{\ell-1} \prod_{i = 0}^{n-1} ({}^{a_i}b_i|_{x.a_n^j})\\
		&\hspace{-0.66em}\begin{array}{llll}
			\begin{array}{l} = \\ \phantom{b} \\ \phantom{\vdots} \\ \phantom{b} \end{array}&
			\begin{array}{rl} b_0\hspace{-0.8em}&|_{x.a_0} \\ \cdot\;\; b_0\hspace{-0.8em}&|_{x.a_na_0} \\ &\vdots \\ \cdot\;\; b_0\hspace{-0.8em}&|_{x.a_n^{\ell-1}a_0} \end{array}&
			\hspace{-1em}\begin{array}{rl} \cdot\;\; b_1\hspace{-0.8em}&|_{x.a_1} \\ \cdot\;\; b_1\hspace{-0.8em}&|_{x.a_na_1} \\ &\vdots \\ \cdot\;\; b_1\hspace{-0.8em}&|_{x.a_n^{\ell-1}a_1} \end{array}&
			\begin{array}{rl} \cdots\phantom{b_0} \\ \cdots\phantom{b_0} \\ \phantom{\vdots} \\ \cdots\phantom{b_0} \end{array}
			\begin{array}{rl} \cdot\;\; b_{n-1}\hspace{-0.8em}&|_{x.a_{n-1}} \\ \cdot\;\; b_{n-1}\hspace{-0.8em}&|_{x.a_na_{n-1}} \\ &\vdots \\ \cdot\;\; b_{n-1}\hspace{-0.8em}&|_{x.a_n^{\ell-1}a_{n-1}}. \end{array}
		\end{array}
	\end{align*}
	Notice that in every column, we are taking the sections of a fixed element $b_i \in B$ at the vertices from the shifted orbit $\orb_{a_n}(x).a_i$. Consequently all sections are taken at different vertices, and there is at most one directed element in every column of the product. Since there are only $n$ columns, by $\syl(g) = \syl(g\|_x)$ every column must contain precisely one directed element, so for all $i$ we have $0 \in \orb_{a_n}(x).a_i$, hence there is some integer $j(i)$ such that $x.a_n^{-j(i)}a_i = 0$. Thus $0.a_i^{-1}a_n^{j(i)} = x$, and ${}^{a_i^{-1}}a_n = {}^{a_i^{-1}a_n^{j(i)}}a_n \in \mathfrak C(a_n, x)$. Lastly we notice that $x.a_n^{j}a_i = 0.{}^{a_i^{-1}}(a_n^{j+j(i)})$, thus
	\[
		g\|_x = \prod_{j = 0}^{\ell-1} \prod_{i = 0}^{n-1} b_i|_{x.a_n^ja_i} = \prod_{j = 0}^{\ell-1} \prod_{i = 0}^{n-1} b_i|_{0.^{a_i^{-1}}(a_n^{j+j(i)})}.\qedhere
	\]
\end{proof}

\begin{proof}[Proof of \cref{Theorem}{thm:A}]
	We establish condition (\ref{eq:lr}) of \cref{Proposition}{prop:torsion_by_length}, for $T = B$.

	Let $x \in X$ and let $g \in G$ be of syllable length $n > 0$ and let
	\begin{equation*}\label{eq:repr g}
		g = {}^{a_0}b_0 \dots {}^{a_{n-1}}b_{n-1}a_n\tag{$\mathsection$}
	\end{equation*}
	such that $n = \syl(g\|_x) = \syl(g)$. If such an element does not exist outside of $B$, (\ref{eq:lr}) is established.
	
	By \cref{Lemma}{lem:non-reduction in one step} the elements ${}^{a_i^{-1}}a_n$ are members of $\mathfrak C(a_n, x)$ for $i \in \{0, \dots, n-1\}$. Now $\orb_{a_n}(x).a_i = \orb_{{}^{a_i^{-1}}a_n}(0)$, hence we write $\ell \defeq \ell_g(x) = \ell_{a_n}(x) = \ell_{{}^{a_i^{-1}}a_n}(0)$. Using the notation of \cref{Lemma}{lem:non-reduction in one step} we have
	\begin{equation*}\label{eq:repr g x}
		g\|_x = \prod_{j = 0}^{\ell-1} \prod_{i = 0}^{n-1} b_i|_{0.^{a_i^{-1}}(a_n^{j+j(i)})} 
		\equiv_{\St_G(1)} \prod_{j = 0}^{\ell-1} \prod_{\substack{i = 0\\j+j(i) \neq \ell}}^{n-1} b_i|_{0.(^{a_i^{-1}}a_n)^{j+j(i)}}.\tag{$\mathsection\mathsection$}
	\end{equation*}
	Since $G$ is strongly orbitwise-abelian we may reorder the product, calculating:
	\begin{align*}
		g\|_x &\equiv_{\St_G(1)} \prod_{i=0}^{n-1} \prod_{\substack{j = 0\\j+j(i) \neq \ell}}^{\ell-1}  b_i|_{0.(^{a_i^{-1}}a_n)^{j+j(i)}}\\
		&= \prod_{i = 0}^{n-1} \prod_{j = 1}^{\ell-1}  b_i|_{0.(^{a_i^{-1}}a_n)^{j}}&& \text{(shifting by $j(i)$)}\\
		&= \prod_{i = 0}^{n-1} \lambda_{b_i}({}^{a_i^{-1}}a_n) && \text{(since $\ell = \ell_{{}^{a_i^{-1}}a_n}(0)$)}\\
		&= \prod_{i = 0}^{n-1} \lambda_{b_i}(a_n, x) && \text{(since $G$ is stable and ${}^{a_i^{-1}}a_n \in \mathfrak C(a_n, x)$)}\\
		&= \lambda_{\prod_{i = 0}^{n-1} b_i}(a_n, x). && \text{(by \cref{Lemma}{lem:orbitwise_abelian_dyn})}
	\end{align*}
	We write $b \defeq \prod_{i = 0}^{n-1} b_i$ to shorten the notation. The $B$-symbols in the middle term of (\ref{eq:repr g x}) representing $g\|_x$ are precisely the $B$-symbols occuring in the word (\ref{eq:repr g}) representing $g$, thus there exists some permutation $\sigma$ of $\{0, \dots, n-1\}$ and some $\hat a_i \in A$ for $i \in \{0, \dots, n-1\}$ such that
	\[
		g\|_x = {}^{\hat a_0}b_{0.\sigma} \dots {}^{\hat a_{n-1}} b_{(n-1).\sigma}\; \lambda_{b}(a_n, x).
	\]
	Since $B$ is abelian, $\prod_{i = 0}^{n-1} b_{i.\sigma} = b$ and we may iterate our calculation for elements $v = x_0 \star \dots \star x_{k-1} \in X^k$, $k \in \N$, such that $\syl(g\|_v) = \syl(g)$. Then
	\[
		g\|_v \bmod{\St(1)} \in \lambda_{b}(\dots (\lambda_{b}(\lambda_{b}(a_n, x_0), x_1), \dots), x_{k-1}) \subseteq \Lambda^k_{b}(a_n).
	\]
	Since $\Lambda_b$ is eventually trivial we conclude that there exist some $k \in \N$ such that for all $v \in X^k$ either $\syl(g\|_v) < \syl(g)$ or $g\|_{v} \in \St_G(1)$. To obtain (\ref{eq:lr}) it is enought to consider all $v$ such that the second case holds. For these, by \cref[(1)]{Lemma}{lem:length_reduction_in_GGS} either $\syl(g\|_{v \star x}) < \syl(g)$ or $g\|_{v \star x} \in B$ holds, since
	\(
		g\|_{v \star x} = g\|_v|_x.
	\)
\end{proof}


\subsection{Discussion of \cref{Theorem}{thm:A}} 
\label{sub:discussion_of_cref_theorem_thm_abelian_case}

Clearly, all \csp groups with abelian rooted group are strongly orbitwise-abelian. Furthermore, if $A$ is abelian, it must act regularly on $X$, hence $\mathfrak C(a, x) = \{ {}^xa \} = \{a\}$. Thus $G$ is stable, and $\lambda_b(a, x) = \lambda_b(a)$ for all $a \in A, b \in B$ and $x \in X$. We record the statement of \cref{Theorem}{thm:A} in this case as a corollary, which turns out to be a version of a result of Bartholdi~\cite{Bar00}*{Theorem 13.4}, with the additional feature that non-cyclic directed groups are allowed. We remark that~\cite{Bar00}*{Theorem 13.4} does not explicitly require that the rooted group is abelian; this appears to be an oversight, see \cref{Example}{eg:non-abelian failure}. Using a similar construction to the one in~\cite{Bar00}, we are able to prove that the condition is not only sufficient but necessary.

\begin{cor}\label{cor:abelian_case}
	Let $G = \la A \cup B \ra$ be a \csp group with abelian and periodic rooted group $A$, and such that either (i) $A$ is of finite exponent, or (ii) $B$ has finite support. Then $G$ is periodic if and only if the dynamical system $(A, \lambda_b)$ is eventually trivial for all $b \in B$.
\end{cor}

\begin{proof}
	Since $A$ is abelian, $B$ is abelian as well. Assuming that (i) or (ii) holds, $B$ is also periodic. Now, identifying the singletons of $\mathcal P(A)$ with $A$, we have $\Lambda_b|_A = \lambda_b$, and the `if'-direction follows from \cref{Theorem}{thm:A}.
	
	For the other implication, assume that there are $a \in A$ and $b \in B$ such that $\lambda_b^n(a) \neq 1_A$ for all $n \in \N$. Consider the orbit length of elements of the form $0^{\ast n}$ under $ba$. From $\ell_{ba}(0^{\ast n}) = \ell_{ba}(0) \cdot \ell_{(ba)\|_0}(0^{\star n-1})$ we conclude that either the orbit lengths are unbounded and $ba$ has infinite order, or that there is some $m \in \N$ such that $\ell_{(ba)\|_{0^{\star n}}}(0) = 1$ for $n > m$. But
	\[
		(ba)\|_0 = b \lambda_b(a),
	\]
	hence $\ell_{(ba)\|_{0^{\ast n}}}(0) = \ell_{\lambda^n_b(a)}(0) > 1$ for all $n \in \N$. Thus, $ba$ has infinite order.
\end{proof}

\cref{Corollary}{cor:abelian_case} may be reformulated to say: A \csp group with abelian rooted group satisfying either (i) or (ii) is periodic if and only if the elements $ba \in B\cdot A$ are of finite order. Motivated by this characterisation, we ask if similar bounds exist for nilpotent and for soluble groups:

\begin{qst}
	For any group $G$, let $I(G)$ be the set of all elements of infinite order. If $G$ is a \csp group, write $B_G(n)$ for the set of elements of syllable length at most $n$. Let $\mathcal C(c)$ be the set of all non-periodic (isomorphism classes of) \csp groups $G = \la A \cup B \ra$ with rooted group nilpotent of class $c$ and satisfying (i) or (ii). Define a function $f_{\mathrm{nil}}: \N \to \N \cup \{\infty\}$ by
	\[
		f_{\mathrm{nil}}(c) = \max\limits_{G \in \mathcal C(c)} \min \{ n \in \N \mid I(G) \cap B_G(n) \neq \emptyset \}.
	\]
	Is $f_{\mathrm{nil}}(c) < \infty$ for all $c \in \N$? One can define an analogous function $_{\mathrm{sol}}$ for soluble groups. Is $f_{\mathrm{sol}}(l) < \infty$ for all $l \in \N$?
\end{qst}

\cref{Theorem}{thm:A} says that $f_{\mathrm{nil}}(1) = f_{\mathrm{sol}}(1) = 1$. We now provide an example of a non-periodic \csp group where the dynamical systems defined by $\lambda_b$ are all eventually trivial, i.e.\ all elements of syllable length less than two are of finite order. The rooted group is isomorphic to the alternating group on four elements acting naturally, implying $f_{\mathrm{sol}}(2) \geq 2$.

\begin{eg}[\cref{Corollary}{cor:abelian_case} cannot be extended to \csp groups with non-abelian rooted group]\label{eg:non-abelian failure}
	Define $s_1 = (0 \; 3 \; 2), s_2 = (0 \; 1 \; 3), s_3 = s_2^{-1}$ and $A = \la s_1, s_2, s_3 \ra \cong \mathrm{A}_4$. Set $b = (b, s_1, s_2, s_3)$ and $B = \la b \ra$. Clearly $b$ has order $3$. The maps $\lambda_b$ and $\lambda_{b^2}$ are shown in \cref{Figure}{fig:lambda}, and they are eventually trivial. Clearly (i) and (ii) are satisfied.

	We now produce an element of infinite order. Let
	\[
		g = s_3s_1bs_3b = (0 \; 2)(1 \; 3) b (0 \; 3 \; 1) b.
	\]
	Then
	\begin{align*}
		g\|_0 &= g^3|_0 = b|_2 b|_2 b|_0 b|_3 b|_1 b|_0 = s_3 b s_3s_1 b\\
		g\|_{0\star 0} &= (s_3 b s_3s_1 b)^3|_0 = b|_3 b|_1 b|_0 b|_2 b|_2 b|_0 = g.
	\end{align*}
	Thus $g$ is not of finite order.
\end{eg}

\begin{figure}
	\centering
	\begin{tikzpicture}[>=stealth,semithick,->,font=\footnotesize, rotate=90]
		\node[state]	(tr)	at (3,-0.5)	{$(\;)$};
		\node[state]	(023)	at (3,-2.5)	{$(0 \; 2 \; 3)$};
		\node[state]	(032)	at (3,1.5)	{$(0 \; 3 \; 2)$};

		\node[state]	(132)	at (5,-2.5)	{$(1 \; 3 \; 2)$};
		\node[state]	(123)	at (5,-0.5)	{$(1 \; 2 \; 3)$};
		\node[state]	(0123)	at (5,1.5)	{$(0 \; 1)(2 \; 3)$};

		\node[state]	(021)	at (7,-0.5)	{$(0 \; 2 \; 1)$};
		\node[state]	(012)	at (7,-2.5)	{$(0 \; 1 \; 2)$};

		\node[state]	(013)	at (3,3.5)	{$(0 \; 1 \; 3)$};
		\node[state]	(0213)	at (5,3.5)	{$(0 \; 2)(1 \; 3)$};
		\node[state]	(031)	at (7,3.5)	{$(0 \; 3 \; 1)$};
		\node[state]	(0312)	at (5,5.5)	{$(0 \; 3)(1 \; 2)$};

		\path	(0312)	edge	(031)
				(031)	edge[bend right]	(0213)
				(0213)	edge[bend right]	(013)
				(013)	edge	(0123)
				(021)	edge[bend right]	(012)
				(012)	edge	(132)
				(0123)	edge	(032)
				(123)	edge	(032)
				(132)	edge	(032)
				(032)	edge[bend right]	(tr)
				(023)	edge[bend right]	(tr);

		\path[dashed,red]
				(0312)	edge	(013)
				(031)	edge	(0123)
				(0213)	edge[bend right]	(031)
				(013)	edge[bend right]	(0213)
				(021)	edge	(123)
				(012)	edge[bend right]	(021)
				(123)	edge	(023)
				(132)	edge	(023)
				(0123)	edge	(023)
				(032)	edge[bend left]	(tr)
				(023)	edge[bend left]	(tr);
	\end{tikzpicture}
	\caption{The maps $\lambda_b$ and $\lambda_{b^2}$ of the group in \cref{Example}{eg:non-abelian failure}. Dashed red arrows represent the latter, solid black arrows the former.}\label{fig:lambda}
\end{figure}
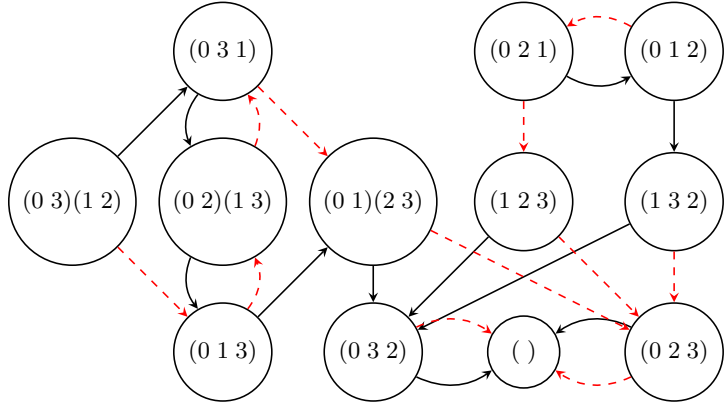

Aside from \csp groups with abelian rooted groups, \cref{Theorem}{thm:A} applies to far more \csp groups. We give an exemplary construction.

\begin{figure}\label{fig:eg thm A}
	\centering
	\begin{tikzpicture}[scale=0.6,x  = {(0.5cm,0.5cm)},
	                    y  = {(0.95cm,-0.25cm)},
	                    z  = {(0cm,0.9cm)}]
		\draw[white] (0,0,-2) -- (0,0,-3);
		\begin{scope}[thick]
			\draw[red] (-1,0,0) -- (1,0,0);
			\draw[green] (0,-1,0) -- (0,1,0);
			\draw[blue] (0,0,-1) -- (0,0,1.25);
		\end{scope}
	
		\fill (0,0,1) circle (.05);
	
		\begin{scope}
			\draw (-1,-1,-1) -- (-1,-1,1);
			\draw (-1,-1,-1) -- (-1,1,-1);
			\draw[dotted] (-1,-1,-1) -- (1,-1,-1);
	
			\draw (-1,-1,1) -- (-1,1,1);
			\draw (-1,-1,1) -- (1,-1,1);
	
			\draw (-1,1,-1) -- (-1,1,1);
			\draw (-1,1,-1) -- (1,1,-1);
	
			\draw[dotted] (1,-1,-1) -- (1,-1,1);
			\draw[dotted] (1,-1,-1) -- (1,1,-1);
	
			\draw (-1,1,1) -- (1,1,1);
			\draw (1,-1,1) -- (1,1,1);
			\draw (1,1,-1) -- (1,1,1);
		\end{scope}
	
		\node at (0,.4,1.3) {$\pi$};	
		\draw[rotate=10] (1.8,-.7) arc (0:320:0.1 and 0.133);	
	\end{tikzpicture}\hspace{1cm}
	\begin{tikzpicture}[scale=0.6,x  = {(0.5cm,0.5cm)},
	                    y  = {(0.95cm,-0.25cm)},
	                    z  = {(0cm,0.9cm)}]
		\draw[white] (0,0,-2) -- (0,0,-3);
		\begin{scope}[thick]
			\draw[teal] (1,-1,-1) -- (-1,1,1);
			\draw[olive] (-1,1,-1) -- (1,-1,1);
			\draw[orange] (1,1,-1) -- (-1,-1,1);
			\draw[magenta] (-1,-1,-1) -- (1.25,1.25,1.25);
		\end{scope}
		\begin{scope}[black!90]
			\draw (-1,-1,-1) -- (-1,-1,1);
			\draw (-1,-1,-1) -- (-1,1,-1);
			\draw[dotted] (-1,-1,-1) -- (1,-1,-1);
	
			\draw (-1,-1,1) -- (-1,1,1);
			\draw (-1,-1,1) -- (1,-1,1);
	
			\draw (-1,1,-1) -- (-1,1,1);
			\draw (-1,1,-1) -- (1,1,-1);
	
			\draw[dotted] (1,-1,-1) -- (1,-1,1);
			\draw[dotted] (1,-1,-1) -- (1,1,-1);
	
			\draw (-1,1,1) -- (1,1,1);
			\draw (1,-1,1) -- (1,1,1);
			\draw (1,1,-1) -- (1,1,1);
		\end{scope}
		
		\node at (1.6,1.5,.6) {$\frac 2 3 \pi$};	
		\draw[rotate=130] (-3.3,1.7) arc (0:320:0.1 and 0.133);
	
	\end{tikzpicture}
	\begin{tikzpicture}[x  = {(0.5cm,0.5cm)},
	                    y  = {(0.95cm,-0.25cm)},
	                    z  = {(0cm,0.9cm)}]
	
		\begin{scope}[scale=2]
		
			\fill[teal!20, thick] (1,0) -- (1,1) -- (0,1) -- (0,0) -- (1,0);
			\fill[orange!20, thick] (1,2) -- (1,3) -- (0,3) -- (0,2) -- (1,2);
			\fill[blue!20, thick] (-1,2) -- (0,2) -- (0,3) -- (-1,3) -- (-1,2);
			\fill[blue!20, thick] (1,2) -- (2,2) -- (2,3) -- (1,3) -- (1,2);
			\fill[red!20, thick] (1,3) -- (1,4) -- (0,4) -- (0,3) -- (1,3);
			\fill[red!20, thick] (1,1) -- (1,2) -- (0,2) -- (0,1) -- (1,1);
		\end{scope}
	
		\begin{scope}[shift={(1,1,0)}, scale=0.3]
			\draw[olive, thick] (-1,1,-1) -- (1,-1,1);
		
			\draw (-1,-1,-1) -- (-1,-1,1);
			\draw (-1,-1,-1) -- (-1,1,-1);
			\draw[dotted] (-1,-1,-1) -- (1,-1,-1);

			\draw (-1,-1,1) -- (-1,1,1);
			\draw (-1,-1,1) -- (1,-1,1);

			\draw (-1,1,-1) -- (-1,1,1);
			\draw (-1,1,-1) -- (1,1,-1);

			\draw[dotted] (1,-1,-1) -- (1,-1,1);
			\draw[dotted] (1,-1,-1) -- (1,1,-1);

			\draw (-1,1,1) -- (1,1,1);
			\draw (1,-1,1) -- (1,1,1);
			\draw (1,1,-1) -- (1,1,1);
		\end{scope}
		\begin{scope}[shift={(1,3,0)}, scale=0.3]
			\draw[blue, thick] (0,0,-1) -- (0,0,1);
		
			\draw (-1,-1,-1) -- (-1,-1,1);
			\draw (-1,-1,-1) -- (-1,1,-1);
			\draw[dotted] (-1,-1,-1) -- (1,-1,-1);

			\draw (-1,-1,1) -- (-1,1,1);
			\draw (-1,-1,1) -- (1,-1,1);

			\draw (-1,1,-1) -- (-1,1,1);
			\draw (-1,1,-1) -- (1,1,-1);

			\draw[dotted] (1,-1,-1) -- (1,-1,1);
			\draw[dotted] (1,-1,-1) -- (1,1,-1);

			\draw (-1,1,1) -- (1,1,1);
			\draw (1,-1,1) -- (1,1,1);
			\draw (1,1,-1) -- (1,1,1);
		\end{scope}
		\node at (1,5,0) {$b$};
		%
		%
		%
		%
		\begin{scope}[shift={(1,7,0)}, scale=0.3]
			\draw[blue, thick] (0,0,-1) -- (0,0,1);
		
			\draw (-1,-1,-1) -- (-1,-1,1);
			\draw (-1,-1,-1) -- (-1,1,-1);
			\draw[dotted] (-1,-1,-1) -- (1,-1,-1);

			\draw (-1,-1,1) -- (-1,1,1);
			\draw (-1,-1,1) -- (1,-1,1);

			\draw (-1,1,-1) -- (-1,1,1);
			\draw (-1,1,-1) -- (1,1,-1);

			\draw[dotted] (1,-1,-1) -- (1,-1,1);
			\draw[dotted] (1,-1,-1) -- (1,1,-1);

			\draw (-1,1,1) -- (1,1,1);
			\draw (1,-1,1) -- (1,1,1);
			\draw (1,1,-1) -- (1,1,1);
		\end{scope}
		\begin{scope}[shift={(3,5,0)}, scale=0.3]
			\draw[blue, thick] (0,0,-1) -- (0,0,1);
		
			\draw (-1,-1,-1) -- (-1,-1,1);
			\draw (-1,-1,-1) -- (-1,1,-1);
			\draw[dotted] (-1,-1,-1) -- (1,-1,-1);

			\draw (-1,-1,1) -- (-1,1,1);
			\draw (-1,-1,1) -- (1,-1,1);

			\draw (-1,1,-1) -- (-1,1,1);
			\draw (-1,1,-1) -- (1,1,-1);

			\draw[dotted] (1,-1,-1) -- (1,-1,1);
			\draw[dotted] (1,-1,-1) -- (1,1,-1);

			\draw (-1,1,1) -- (1,1,1);
			\draw (1,-1,1) -- (1,1,1);
			\draw (1,1,-1) -- (1,1,1);
		\end{scope}
		\begin{scope}[shift={(-1,5,0)}, scale=0.3]
			\draw[blue, thick] (0,0,-1) -- (0,0,1);
		
			\draw (-1,-1,-1) -- (-1,-1,1);
			\draw (-1,-1,-1) -- (-1,1,-1);
			\draw[dotted] (-1,-1,-1) -- (1,-1,-1);

			\draw (-1,-1,1) -- (-1,1,1);
			\draw (-1,-1,1) -- (1,-1,1);

			\draw (-1,1,-1) -- (-1,1,1);
			\draw (-1,1,-1) -- (1,1,-1);

			\draw[dotted] (1,-1,-1) -- (1,-1,1);
			\draw[dotted] (1,-1,-1) -- (1,1,-1);

			\draw (-1,1,1) -- (1,1,1);
			\draw (1,-1,1) -- (1,1,1);
			\draw (1,1,-1) -- (1,1,1);
		\end{scope}
	\end{tikzpicture}
	\caption{The action of $\Alt(4)$ on the cube $X$ and the action of $b$ on $X^*$ described in \cref{Example}{eg:non-abelian strongly orbitwise-abelian stable}.}
\end{figure}
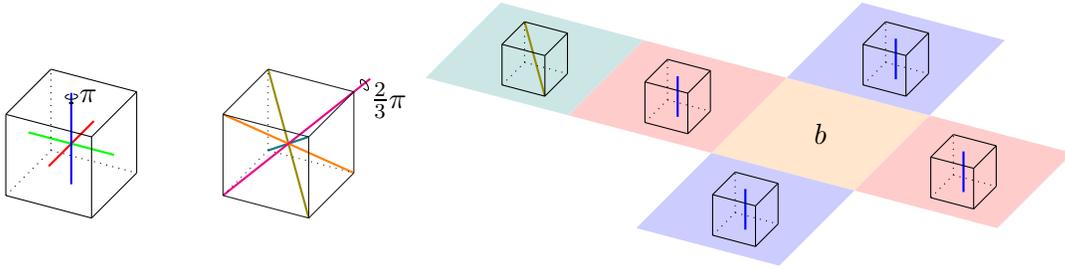

\begin{eg}[A non-abelian strongly orbitwise-abelian and stable \csp group]\label{eg:non-abelian strongly orbitwise-abelian stable}
	Let $\Sym(4)$ be the symmetric group on four letters. The set $X$ of transpositions in $\Sym(4)$ is a conjugacy class of $\Sym(4)$, and the alternating group $A = \Alt(4)$ acts faithfully and transitively on $X$ by conjugation. Equvialently, the alternating group acts on the faces of a cube by rotations in the following way. The double transpositions rotate the cube by $\pi$ along an axis defined by the center points of opposing sides, and the three-cycles rotate it by $\frac 2 3 \pi$ along an axis going through two opposing corners, see \cref{Figure}{fig:eg thm A}.
	
	We fix the transposition $(1\,2)$, resp.\ the top face of the cube, as the distinguished letter. Clearly $\st_{A}(1\,2) = \la (1\,2)(3\,4) \ra$ is of order two. We write $V$ for the Klein four group consisting of the double transpositions. Define an automorphism $b$ by
	\[
		b = ((1\,2): b, \; (3\,4): (1\,2\,3), \; x: (1\,2)(3\,4) \text{ for } x \in X\setminus\{(1\,2),(3\,4)\}).
	\]
	Compare again with \cref{Figure}{fig:eg thm A}.
	The \csp group $G$ defined by $A$ and $B = \la b \ra$ is stable, strongly orbitwise-abelian and the dynamical system defined by $\Lambda_{b'}$ for any $b' \in B$ is eventually trivial, hence periodic by \cref{Theorem}{thm:A}. Our discussion will show that minor variations in the definition of $b$ lead to the same conclusion, e.g.\ one may replace $b|_{(3\,4)}$ by an arbitrary three-cycle, or define the sections of $b$ on the blue (resp.\ the red) faces to be trivial. We obmit more careful calculations that show, again using \cref{Theorem}{thm:A}, that there are further examples of \csp groups with rooted group $A$.
	We first show that $G$ is strongly orbitwise-abelian. Let $a \in V$. Since $V$ is normal, $\mathfrak C(a,x) = {}^{\mathrm{mp}_A( (1\,2) ,x)}a \subset V$, and since it is abelian, it is a singleton set. Now the image of $(1\,2)$ under an element of $V$ is either $(1\,2)$ or $(3\,4)$ (the bottom side of the cube). Consequently
	\[
		\mathfrak X(a,x) \in \{ \{(1\,2)\}, \{(3\,4)\}\}
	\]
	for all $a \in V, x \in X$. Now let $a \in A$ be an element of order three. The set $T$ of such elements is a normal subset, hence $\mathfrak C(a,x) \subseteq T$. Now $T$ acts (as a set) non-transitively on $X$, since $(3\,4) \notin T.(1\,2)$. Viewing at the cube, it is apparent that the bottom and top faces cannot be interchanged by an element of $T$. At the same time, no element of $T$ has any fixed points. Consequently,
	\[
		\mathfrak X(a,x) \in \{ \{(1\,3), (1\,4)\}, \{(2\,3), (2\,4)\}\}
	\]
	for all $a \in T, x \in X$. These are sets of opposing sides. A calculation shows that $\mathfrak X(a,x)$, with $a$ ranging over $A$, takes all values
	\[
		\{\{(1\,2)\}, \{(3\,4)\}, \{(1\,3), (1\,4)\}, \{(2\,3), (2\,4)\}\}.
	\]
	These correspond to the colouring in \cref{Figure}{fig:eg thm A}. Clearly the sections within such a set commute, hence $G$ is orbitwise-abelian.
	
	Since $\mathfrak C(a,x)$ is a singleton for $a \in V$, the group $G$ being stable is equivalent to $\lambda_b$ taking the same values for all $a \in \mathfrak C(t, x)$, $t \in T$, $x \in X$, $b' \in B$. The two elements in any $\mathfrak C(t, x)$ are $\st_A(1\,2)$-conjugate, hence
	\[
		\mathfrak C(t,x) \in \{ \{ (1\,2\,3), (1\,4\,2) \}, \{ (1\,3\,2), (1\,2\,4) \}, \{ (1\,3\,4), (2\,4\,3) \}, \{ (1\,4\,3), (2\,3\,4) \} \}
	\]
	and all these values are taken. Any such pair corresponds to two orbits of $(1\,2)$ going around opposite corners of the top face, i.e.\ a red and a blue face. The product over the sections is hence $((1\,2)(3\,4))^2 = 1_A$ in all cases. Thus $G$ is stable, and furthermore for all $t \in T$
	\[
		\Lambda_{b}(t) = \Lambda_{b^2}(t) = \{1_A\}.
	\]
	Finally, the conjugate of $(1\,2)$ under a double transposition is either $(1\,2)$ or $(3\,4)$. Hence
	\[
		\Lambda_{b}(a) = \{b|_{(1\,2)}, b|_{(3\,4)} \} = \{1_A, (1\,2\,3)\} \quad\text{ and }\quad \Lambda_{b^2}(a) = \{1_A, (1\,3\,2) \}
	\]
	for all $a \in V$, and by $\Lambda_b^2(a) = \Lambda_{b^2}^2(a) = \{1_A\}$ we see that $\Lambda_{b'}$ is eventually trivial for all $b' \in B$. Thus $G$ is periodic by \cref{Theorem}{thm:A}.
\end{eg}


\subsection{Proof and discussion of \cref{Theorem}{thm:B}} 
\label{sub:proof_of_cref_theorem_thm_general_periodicity_preserving_thm}

If we drop the assertion that $G$ is strongly orbitwise-abelian, the argument used in the proof of \cref{Theorem}{thm:A} to reduce the set of candidates of infinite order to elements $g$ with $\syl(g) \leq 1$ fails in the general case, as illustrated by \cref{Example}{eg:non-abelian failure}. But we may still consider the ``coarser'' dynamical system defined by the maps $\Sigma_B$. This gives us a broadly applicable sufficent condition for periodicity, that does not require $B$ to be abelian. Also, one does not need to check the full system $\Sigma_B$, which may be replaced with $\Sigma_S$, for a generating set $S \subseteq B$.

\begin{proof}[Proof of \cref{Theorem}{thm:B}]
	We aim to apply \cref{Proposition}{prop:torsion_by_length} with $T = B$. Let
	$$g = ({}^{a_{0}}b_{0})\dots({}^{a_{n-1}}b_{n-1})a_n  \in G\setminus A$$
	be an element of syllable length $n > 0$ and $x \in X$ such that $\syl(g\|_x) = n$. Write $\ell \defeq \ell_g(x)$. By \cref{Lemma}{lem:non-reduction in one step} we find that ${}^{a_i^{-1}}a_n \in \mathfrak C(a_n,x)$ for all $i \in \{0, \dots, n-1\}$ and there are integers $j(i)$ such that
	\[
		g\|_x = \prod_{j = 0}^{\ell-1} \prod_{i = 0}^{n-1} b_i |_{0.^{a_i^{-1}}a_n^{j+j(i)}}.
	\]
	Decompose all $b_i$ as a product of generators in $S$, writing $b_i = \prod_{m = 0}^{n_i-1} b_{i, i_m}$. Then calculating modulo $\St(1)$ and $$N \defeq \la b|_y \mid y \in \mathfrak X(a_n, x), b \in S \ra' \leq A'$$
	we find
	\[
		g\|_x \equiv_{\la N, \St(1) \ra} 
		\prod_{i = 1}^{n-1} \prod_{m = 0}^{n_i - 1} \left(\prod_{j = 1}^{\ell-1} b_{i,i_m}|_{0.^{a_i^{-1}}a_n^j}\right).
	\]
	Thus $g\|_x |^\epsilon \in \sigma_S(a_n, x)$.

	Let $k \in \N$ be such that $\Sigma_S^k(a_n) = \{1_A\}$, and let $v \in X^k$. If $\syl(g\|_v) = n$, all intermediate stabilised sections must have syllable length $n$ and we obtain $g\|_v|^\epsilon \in \Sigma_S^k(a_n) = \{1_A\}$ by the calculation above. Let $y \in X$ be any letter. Then
	\[
		\syl(g\|_{vy}) = \syl(g\|_v|_y) < \syl(g\|_v) = n,
	\]
	or $g\|_{vy} \in B$ by \cref[(1)]{Lemma}{lem:length_reduction_in_GGS}. Thus, by \cref{Proposition}{prop:torsion_by_length}, $G$ is periodic.
\end{proof}

\begin{rmk}\label{rmk:regular case}
	The special case of a \csp group where the rooted group $A$ acts regularly (on itself with the distinguished letter $1_A$) was first considered by Gupta and Sidki~\cites{GS83,GS84}. They considered automorphisms defined by \textit{decorating functions}, roughly corresponding to directed elements, and gave a list of five conditions for such a group $G = \la A \cup \{\delta\} \ra$ with periodic rooted group $A$ and an automorphism $\delta$ defined by a decorating function (i.e.\ a \csp~group with cyclic directed group) to be periodic, in this case, the decorating function is called \textit{periodicity preserving}. Translated into the language of \csp groups, the five conditions are:
	\begin{enumerate}
		\item The element $\delta$ is directed.
		\item The group $G$ is a \csp group.
		\item The element $\delta$ has finite support.
		\item The group $G$ is orbitwise-abelian.
		\item For all $a \in A \setminus \{ 1_A\}$ we have $\prod_{a' \in \la a \ra \setminus \{ 1_A\}} \delta|_{a'} = 1_A$.
	\end{enumerate}
	Note that by the fourth condition, the product in the fifth condition is well-defined.

	We argue that the conditions of \cref{Theorem}{thm:B} are implied by the five conditions of Gupta and Sidki. Using \cref{Theorem}{thm:B}, we see from (1), (2) and (3) that we only need to prove that $\Sigma_S$ is eventually trivial for some generating set $S$ of $B$. Naturally, we set $S = \{ \delta \}$.

	Since $A$ acts regularly, $\mathfrak C(a, x) = \{{}^xa\}$ and $\mathfrak X(a,x) = \orb_{{}^xa}(0) \setminus \{ 0 \}$. For all $a \in A$ and non-trivial $x \in A$, we thus have
	\[
		\la \prod_{c' \in \la c \ra  \setminus \st_A(0)} \delta|_{c'} \mid c \in \mathfrak C(a, x) \ra = \la \prod_{a' \in \la {}^xa \ra \setminus \{ 1_A\}} \delta|_{a'} \ra \overset{\text{(5)}}{=} \{1_A\},
	\]
	and the group in \cref{Definition}{def:sigma} of $\sigma_S(a, x)$ that is generated by products is trivial. Since $A$ acts regularly, (4) implies that $G$ is strongly orbitwise-abelian, hence the derived subgroup in the definition of $\sigma_S(a, x)$ is trivial as well. Thus the dynamical system defined by
	\[
		\Sigma_S(a) = \bigcup_{x \in X} \sigma_S(a, x) = \bigcup_{x \in X} \{ 1_A \}
	\]
	is eventually (or rather, immediately) trivial. Consequently, the conditions of Gupta and Sidki appear as a special case of \cref{Theorem}{thm:B}.
\end{rmk}

Gupta and Sidki proved that nearly all (precisely all non-dihedral, non-cyclic) finitely generated periodic groups allow a periodicity preserving function. In case of a finite group $A$ generated by a set of non-involutions $N$, one may define such a function $\beta: A \setminus \{ 1_A\} \to A$ by
\[
	\beta(x) = \begin{cases}
		x &\text{ if }x \in S,\\
		1_A &\text{ otherwise.}
	\end{cases}
\]
It is an interesting question which finite permutation groups may occur as the local action (equivalently, as rooted subgroup) of a periodic \csp group, cf.\ \cref{Subsection}{sub:limitations_and_symmetric_groups}. In combination with \cref{Corollary}{thm:basilica groups}, which we will discuss later, this shows that every group of the form $A^{s} \curvearrowleft A^{\wr s}$ for $A$ a non-dihedral, non-cyclic finite group allows a periodic \csp group that features $A^{\wr s}$ (with this action) as its rooted subgroup. It is not true that every transitive permutation group allows a periodic \csp group.

Even if we restrict to the regular case, \cref{Theorem}{thm:B} is very flexible.
As a demonstration, we prove the following corollary.

\begin{cor}
	Let $A$ be a finite periodic group that decomposes as a semi-direct product $A = N \rtimes \la g \ra$ with a cyclic quotient of order greater than two. There exists a finitely generated infinite periodic group generated by two isomorphic copies of $A$.
\end{cor}

\begin{proof}
	Assume that $A = N \rtimes \la g \ra$ and let $S \subseteq A$ be a generating system containing $g$ such that $S\setminus\{g\}$ generates $N$. We construct a group $B \in \Aut(A^*)$ isomorphic to $A$ by defining 
	\[
		b_s|_g = s, \quad b_s|_{g^{-1}} = \begin{cases}
			1_A &\text{ if } s \neq g,\\
			g^{-1} &\text{ if } s = g,
		\end{cases}, \quad b_s|_x = 1_A, \quad\text{ for }x \in A \setminus\{1_A, g, g^{-1}\},
	\]
	for every $s \in S$. Set $B = \la b_s \mid s \in S \ra$. This group is isomorphic to the group generated by
	\[
		\{ (s, 1_A) \mid s \in S\setminus\{g\} \} \cup \{ (g, g^{-1}) \}.
	\]
	But since $N = \la S\setminus\{g\} \ra$, this group is isomorphic to $A$.
	We now prove that the \csp group defined by $A$ and $B$, with $A$ acting regularly on itself, is periodic. As $B$ is finite, by \cref{Theorem}{thm:B} it is sufficient to prove that the dynamical system defined by $\Sigma_{S_B}$ is eventually trivial, where we set $S_B = \{ b_s \mid s \in S\}$.
	
	Let $x,y \in A$. We calculate $\sigma_{S_B}(x,y)$. First observe that $\mathfrak C(x,y) = \{ {}^yx \}$, since $A$ acts regularly. There are two cases, depending if $a \in \la {}^y x \ra$. In case $g \in \la {}^yx \ra$, also $g^{-1}$ is in $\la {}^y x \ra$, and the products generating the subgroup on the left in the equation (\ref{eg:sigma}) defining $\sigma_{S_B}(x,y)$ in \cref{Definition}{def:sigma} are $s \in S\setminus\{g\}$ and $1_A$, while the generators of the subgroup on the right are the members of $S$. Then, 
	\[
		\sigma_{S_B}(x,y) = N \cdot A' = N.
	\]
	In case $g \notin \la {}^yx \ra$, all sections $b_s|_{z}$ for $z \in \la {}^yx \ra \setminus \st_A(0)$ and $s \in S$ are trivial, hence
	\[
		\sigma_{S_B}(x, y) = \{1_A\}.
	\]
	Thus for any $x \in A$, we have $\Sigma_{S_B}(x) \leq N$.
	
	Now assume $x \in N$. Since $N$ is normal, all conjugates ${}^yx$ are again in $N$, and $g \notin \la {}^yx \ra$. Thus
	\[
		\Sigma_{S_B}(x) = \{1_A\},
	\]
	and $\Sigma_{S_B}$ is eventually trivial.
\end{proof}

Using \cref{Theorem}{thm:B}, it is also not too difficult to find examples of periodic \csp groups where the action of the rooted group is not regular, a previously unrecorded phenomenon.

\begin{eg}\label{eg:B}
	We define a family of \csp groups satisfying the condition of \cref{Theorem}{thm:B}. Let $p$ be an odd prime and identify $X = \{ 0, \dots, 2p-1\}$ with the vertices of the regular $2p$-gon. We fix the rotation $r = (0 \, 1 \, \dots \, 2p-1)$ and the reflexion $s = (1 \, 2p-1)( 2 \, 2p-2) \dots (p-1 \, p + 1)$ in the axis through $0$ and $p$ as the generators of the rooted group $A \cong D_{2p}$, the dihedral group of order $4p$.

	Define an automorphism $b$ by
	\[
		b|_x = \left\{\begin{array}{lllll}
			r &\text{ if }x \equiv_2 1 \text{ and } x \neq p, && 1_A &\text{ if } x = p,\\
			s &\text{ if }x \equiv_2 0 \text{ and } x \neq 0, && b &\text{ if } x = 0.
		\end{array}\right.
	\]
	See \cref{Figure}{fig:hex} for an illustration of this automorphism. Clearly, the first layer sections of $b$ generates $A$. By considering seperately the elements of order $2$, $p$ and $2p$ in $A$, we show that the \csp group defined by $A$ and $\la b \ra$ gives rise to an eventually trivial dynamical system $\Sigma_{\{b\}}$ and is thus periodic.

	Let $a \in A$ be an element of order $p$, i.e.\ $a \in \la r^2 \ra$. The orbit of $0$ under any such element is the set of even elements in $X$. Excluding $0$, all sections at even elements are equal to $s$, hence for all $x \in X$
	\[
		\sigma_{\{b\}}(a, x) = \la s^{p-1} \ra \cdot \la s \ra' = \{1_A\}.
	\]
	Let $a \in A$ be an element of order $2p$, i.e.\ of the form $r^i$ for some odd $i$ that is not a multiple of $p$. The orbit of $0$ under any conjugate of $a$ is the full set $X$, and we find
	\[
		\sigma_{\{b\}}(a, x) = \la r^{p-1} \cdot s^{p-1} \ra \cdot \la s, r \ra' = \la r^2 \ra
	\]
	for any $x \in X$, hence $\Sigma_{\{b\}}(a) = \la r^2 \ra$, which becomes trivial after another step as we saw above.

	Let $a \in A$ be an element of order $2$. Then either $a = sr^i$ for some $i \not\equiv_{2p}0$, $a = s$ or $a = r^p$. In the third case, the element $r^p$ is central, hence $\mathfrak C(r^p, x) = \{ r^p \}$ for all $x \in X$ and $\mathfrak X(r^p,x) = \{p\}$, thus $\Sigma_{\{b\}}(r^p) = \sigma_{\{b\}}(r^p,x) = \{ 1_A \}$. In the second case, the element $s$ fixes the distinguished letter and we also find $\Sigma_{\{b\}}(s) = \{1_A\}$. In the first case, $\mathfrak C(sr^i, x) = \{sr^i, sr^{i-2} \}$ and
	\[
		\sigma_{\{b\}}(sr^{i}, x) = \begin{cases}
			\la r \ra \cdot \la r \ra' = \la r \ra 		&\text{ if } i \equiv_2 1,\\
			\la s \ra \cdot \la s \ra' = \la s \ra		&\text{ if } i \equiv_2 0,
		\end{cases}
	\]
	since $0.sr^i = i \equiv_2 0$ if and only if $0.sr^{-i} = 2p-i \equiv_2 0$. Thus the system $\Sigma_{\{b\}}$ is eventually trivial.
\end{eg}

Obviously there are many variants of this construction, e.g.\ replacing $s$ with another involution of $A$.

Furthermore, using similar arguments, periodic \csp groups with rooted group isomorphic to dihedral groups of even exponent acting naturally may easily be found.

\begin{qst}
	 As described by Grigorchuk in~\cite{Gri00}, a \GGS group on a $p$-adic rooted regular tree is periodic if and only if its defining vector $\underline e$ satisfies the equation $\sum_{i = 1}^{p-1} e_i \equiv_p 0$. Hence for \csp groups with rooted group $C_p$, there is a criterion for periodicity that can be stated rather easily.
	 
	 Is there a similar characterisation for the periodic groups among the \csp groups with rooted group $D_{2p}$?
\end{qst}

\begin{figure}
	\centering
	\begin{tikzpicture}[scale=1.5]
		\fill[teal!20] (0:.5) \foreach \x in {60,120,...,360} { -- (\x:.5)};
		
		\foreach \x in {0,1,...,5} {
			\draw[black!60] (\x*60:.5) -- (\x*60:.7);
			\node at (\x*60:.33) {{\x}};
			\fill (\x*60:.5) circle (.03);
		}
		
		\fill[shift={(0:1)}, orange!20] (0:.3) \foreach \y in {60,120,...,360} { -- (\y:.3)};
		\node at (0:1) {$b$};
		
		\fill[shift={(60:1)}, red!20] (0:.3) \foreach \y in {60,120,...,360} { -- (\y:.3)};
		\draw[shift={(0:.4)}, ->] (60:1) arc (0:60:.3);
		\node at (60:1) {$r$};
		
		\fill[shift={(120:1)}, blue!20] (0:.3) \foreach \y in {60,120,...,360} { -- (\y:.3)};
		\draw[shift={(120:1)}] (0:.4) -- (180:.4);
		\node[shift={(90:.2)}] at (120:1) {$s$};
		
		\fill[shift={(180:1)}, teal!20] (0:.3) \foreach \y in {60,120,...,360} { -- (\y:.3)};
		\node at (180:1) {$1_A$};
		
		\fill[shift={(240:1)}, blue!20] (0:.3) \foreach \y in {60,120,...,360} { -- (\y:.3)};
		\draw[shift={(240:1)}] (0:.4) -- (180:.4);
		\node[shift={(90:-.2)}] at (240:1) {$s$};
				
		\fill[shift={(300:1)}, red!20] (0:.3) \foreach \y in {60,120,...,360} { -- (\y:.3)};
		\draw[shift={(0:.4)}, ->] (300:1) arc (0:60:.3);
		\node at (300:1) {$r$};
	\end{tikzpicture}\hspace{1cm}
	\begin{tikzpicture}[>=stealth,semithick,->,font=\footnotesize]		
		\node	(1A)	at (0,0)	{$\{1_A\}$};
		
		\node	(r2)	at (-2,1)	{$\{r^2\}$};
		\node	(r4)	at (-1,1)	{$\{r^4\}$};
		\node	(r2g)	at  (0,1)	{$\la r^2 \ra$};
		\node	(r3)	at  (1,1)	{$\{r^3\}$};
		\node	(s)		at  (2,1)	{$\{s\}$};
		\node	(sg)	at  (3,1)	{$\la s \ra$};
		
		\node	(r)		at (-1,2)	{$\{r\}$};
		\node	(rg)	at  (0,2)	{$\la r \ra$};
		\node	(r5)	at  (1,2)	{$\{r^5\}$};
		\node	(sr2) 	at  (2,2)	{$\{sr^2\}$};
		\node	(sr4)	at  (3,2)	{$\{sr^4\}$};

		\node	(sr)	at (-1,3)	{$\{sr\}$};
		\node	(sr3)	at  (0,3)	{$\{sr^3\}$};
		\node	(sr5)	at  (1,3)	{$\{sr^5\}$};

		\path	(r2)	edge	(1A)
				(r4)	edge	(1A)
				(r2g)	edge	(1A)
				(r3)	edge	(1A)
				(s)		edge	(1A)
				(sg)	edge[bend left]	(1A);
				
		\path	(r)		edge	(r2g)
				(rg)	edge	(r2g)
				(r5)	edge	(r2g)
				(sr2) 	edge	(sg)
				(sr4)	edge	(sg);
		
		\path	(sr)	edge	(rg)
				(sr3)	edge	(rg)
				(sr5)	edge	(rg);
		\end{tikzpicture}
	\caption{The generator $b$ of \cref{Example}{eg:B} for $p = 3$ and the part of the dynamical system defined by $\Sigma_b$ reached by singletons.}\label{fig:hex}
\end{figure}
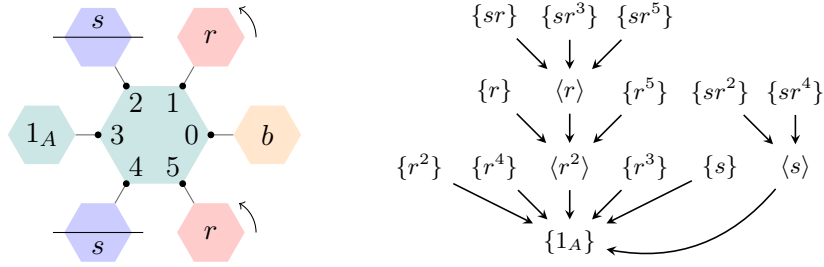



\section{Applications and examples} 
\label{sec:applications_and_examples}

In the previous section we have seen that for nearly all finitely generated periodic groups acting regulary on themselves, we find directed groups such that the resulting \csp group is periodic. Now we will give an example how to construct \csp groups acting locally non-regularly, that fulfill the conditions of \cref{Theorem}{thm:B}, and at the same time explore its limitations.

\subsection{Rooted groups where \cref{Theorem}{thm:A} and \ref{thm:B} cannot be applied} 
\label{sub:limitations_and_symmetric_groups}

We make the following observation.
Let $G$ be a \csp group with rooted group $A \leq \Sym(X)$ that is perfect and contains an element $a \in A$ acting transitively on $X$ (i.e.\ an $|X|$-cycle). Clearly $G$ is not orbitwise-abelian, so \cref{Theorem}{thm:A} cannot be applied to $G$. Furthermore, let $S \subseteq B$ be any generating set for $B$. Then the sections $b|_x$ for $x \in \dot X$ and $b \in S$ generate $A$, and since $\sigma_S(a, 0)$ contains the derived subgroup of the group generated by these sections, we have $\sigma_S^k(a, 0) = A$ for all $k > 0$. Consequently, $G$ does not fulfill the conditions of \cref{Theorem}{thm:B}.

In view of this limitation, we ask:

\begin{qst}
	Are there any \csp groups based on a rooted group $A \cong \Sym(X)$ acting naturally on $X$ that are periodic?
\end{qst}

In regards to this question we point out that it is possible to exclude certain special cases. It is well known that the only \csp group on a two element set is infinite dihedral. An extensive brute-force computer calculation conducted by the author using GAP~\cite{GAP} show that there is no \csp group with rooted group $\Sym(3)$ acting naturally.

We consider another special case, and write $X = \{0, \dots, m-1\}$. Assume that there is a generating set $S$ of $B$ such that the collection of sections $b|_x$ for $b \in S$ and $x \in \dot X$ contains the cycle $a_m = (0 \, 1 \,\dots\, m-1)$ exactly once, the transposition $a_2 = (0 \; 1)$ exactly twice, and nothing more. Then $G$ cannot be periodic:

Let $b \in S$ be the unique element with a section equal to $a_m$. There are three cases: Either $b$ has (aside from the sections $a_m$ and $b$) no, one, or two sections equal to $a_2$. If there are none, $a_mb$ is of infinite order, since $\lambda_b(a_m) = a_m$. If there are two, the element $\lambda_b(a_m)$ is a product of $a_2, a_2$ and $a_m$, hence either equal to $a_m$ or to ${}^{a_2}a_m$. In the first subcase $a_mb$ is of infinite order. In the second subcase, $\lambda_{b}(^{a_2}a_m)$ is again equal to $a_m$ or ${}^{a_2}a_m$. Thus either $\lambda_b$ or $\lambda_b^2$ have a fixed point, and $a_mb$ is of infinite order. Lastly, if there is precisely one section of $b$ equal to $a_2$, there is another generator $b' \in S$ with exactly one section equal to $a_2$. Now $\lambda_{bb'}(a_m)$ is either $a_m$ or ${}^{a_2}a_m$, and we may argue as before.

Clearly, similar methods can be applied to exclude other generating sets of $\Sym(X)$, which feature sparse section decompositions.


\subsection{Basilica groups of \GGS groups} 
\label{sub:basilica_groups_of_ggs_groups}

We now turn our attention to the \textit{$s$\textsuperscript{th} Basilica groups} of \csp groups resulting from the Basilica construction introduced in~\cite{PR21}. To every group of tree automorphisms $H$ one can associate a family of Basilica groups, all sharing certain properties with $H$. Most famously, the second Basilica group of the dyadic odometer is the (classical) Basilica group defined by Grigorchuk and \.{Z}uk~\cite{GZ02}, while the Basilica groups of spinal groups are again spinal, altough on another tree; cf.\ \cite{PR21}*{Proposition 3.9}. We prove that if a \csp group $G$ satisfies the conditions of \cref{Theorem}{thm:B}, its Basilica groups do so as well, providing many examples of periodic \csp groups with non-regular rooted action. In particular, this is also of interest in the context of Basilica groups, as it is not known if Basilica groups of periodic groups are again periodic.

Instead of recalling the general definition of the $s$\textsuperscript{th} Basilica group, we provide an ad hoc description for \csp groups.

\begin{definition}
	Let $G = \la A \cup B \ra$ be a \csp group acting on the tree $Y$ and $s \in \N$. For $i \in \{0, \dots, s-1\}$ and for every $a \in A$ define an element $a_i \in A^{\wr s}$ by the following action on $X = Y^s$:
	\[
		(y_0 \star\dots \star y_{s-1}).a_i \defeq y_0 \star\dots \star y_{i-1} \star y_i.a \star y_{i+1} \star \dots \star y_{s-1}.
	\]
	Then $\la a_i \mid a \in A, i \in \{0, \dots, s-1\} \ra \cong A^{\wr s}$ acts transitively on $X$. For every $i \in \{ 0, \dots, s-1\}$ and $b \in B$ define an automorphism of $X^*$ by
	\[
		b_i|_u = \begin{cases}
			a_{i} &\text{ if }u = 0^{\star(i)}\star y \star 0^{\star(s-1-i)} \text{ and } a = b|_y,\\
			\id &\text{ otherwise.}
		\end{cases}
	\]
	Defining $A_i = \{ a_i \mid a \in A\}$ and $B_i$ equivalently, we define the $s$\textit{\textsuperscript{th} Basilica group} $\operatorname{Bas}_s(G)$ as the \csp group on the tree $X^*$ defined by the rooted group $\la \, \bigcup_{i = 0}^{s-1} A_i \,\ra$ and the directed group $\la \, \bigcup_{i = 0}^{s-1} B_i \,\ra$.
\end{definition}

It is a general fact that the directed group of a Basilica group of a \csp group (or, more generally, any spinal group) $G$ is the $s$-fold direct product of the directed group of $G$, cf.\ \cite{PR21}*{Proposition 3.9}. Let $T \subset G$ be a generating set of $G$. Then the set
\[
	\mathbb T = \{ t_i \mid i \in \{ 0, \dots, s-1 \}, t \in T \}
\]
is a generating set for $\operatorname{Bas}_s(G)$.

Note that $\operatorname{Bas}_s(G)$ is a \csp group on the tree $X^*$, where the alphabet $X$ itself is the finite rooted tree $Y^s$. To distinguish the two tree-structures we omit the symbol $\star$ for elements in the finite tree from here on.

\begin{cor}\label{thm:basilica groups}
	Let $G$ be a \csp group satisfying the conditions of \cref{Theorem}{thm:B} for the generating system $T \subseteq G$. Then $\operatorname{Bas}_s(G)$ is periodic for all $s \in \N$.
\end{cor}

\begin{proof}
	Since the rooted group of $\operatorname{Bas}_S(G)$ is isomorphic to $A^{\wr s}$ and the directed group of $\operatorname{Bas}_S(G)$ is isomorphic to $B^s$, conditions (i) and (ii) transfer to $B^s$, and by \cref{Theorem}{thm:B} it is enough to prove that $\Sigma_{\mathbb T}$ is eventually trivial.

	The rooted group $A^{\wr s}$ acts on the finite rooted regular tree $X = Y^s$. Write $S(i)$ for the $i$\textsuperscript{th} level stabiliser of $A^{\wr s}$, with is isomorphic to $(A^{\wr s-i})^{|Y|^i}$, and set $S(s) = \{1_{A^{\wr s}}\}$. Clearly $A_j \subseteq S(i)$ if and only if $j \geq i$.

	Let $k \in \N$ be such that $\Sigma_T^k(A) = \{1_A\}$. We shall prove that $\Sigma^{s\cdot k}_{\mathbb T}(A^{\wr s})$ is trivial by proving that
	\[
		\Sigma_{\mathbb T}^k(S(i)) \subseteq S(i+1)
	\]
	for all $i \in \{0, \dots, s-1\}$. Thus let $g \in S(i)$ and let $x = y_0 \dots y_{s-1} \in Y^s$. To calculate 
	\[
		\sigma_{\mathbb T}(g, x) = \la \prod_{c' \in \la c \ra  \setminus \st_{A^{\wr s}}(0^s)} t|_{0^{s}.c'} \mid c \in \mathfrak C(g, x), t \in \mathbb T \ra \cdot \la t|_{x'} \mid x' \in \mathfrak X(g, x), t \in \mathbb T \ra'.
	\]
	we first notice that $\mathfrak C(g, x)$ is a subset of the conjugates of $g$, hence since $S(i)$ is normal, $\mathfrak C(g, x) \subseteq S(i)$. Thus for any $c \in \mathfrak C(g, x)$ the orbit of the distinguished letter is of the form
	\[
		\orb_c(0^{s}) \subseteq 0^{i} Y^{s-i},
	\]
	and by the definition of $t_j$ for $t \in T$, $j \in \{0, \dots, s-1\}$, the sections of $t_j$ along $\orb_c(0^{s})\setminus\{ 0^s\}$ are trivial if $j < i$. Now $B_j|_{x'} \subseteq B_j \cup A_j \subseteq S(i)$ for all $x' \in X$ implies $\sigma_{\mathbb T}(g, x) \subseteq S(i)$.

	We now calculate modulo $S(i+1)$. Consequently we may ignore all sections of elements $t_j$ with $t \in T$ and $j \neq i$, and get
	\begin{align*}
		\sigma_{\mathbb T}(g, x) &\equiv_{S(i+1)} \la \prod t_i|_{0^{s}.c} \mid c \in \mathfrak C(g, x), t \in T \ra \cdot \la t_i|_{x'} \mid x' \in \mathfrak X(g, x), t \in T \ra'\\
		&\equiv_{S(i+1)} \sigma_T(g \bmod{S(i+1)}, y_0\dots y_i).
	\end{align*}
	Thus $\Sigma_{\mathbb T}^k(S(i)) \subseteq S(i+1)$, and $\Sigma_{\mathbb T}$ is eventually trivial.
\end{proof}


\subsection{A periodic \csp group satisfying neither \cref{Theorem}{thm:A} nor \ref{thm:B}} 
\label{sub:an_example}

At last, we demonstrate that, in contrast to \cref{Corollary}{cor:abelian_case} for abelian rooted groups, neither \cref{Theorem}{thm:A} nor \cref{Theorem}{thm:B} nor the union of their scopes provide necessary conditions for periodicity. Indeed, we construct a periodic \csp group that is subject to neither of the sets of conditions. 

Let $X = \{0, 1, 2, 3 \}$ be the alphabet and let $A \leq \Sym(X)$ be generated by $s = (1 \, 3)$ and $r = (0 \, 1 \, 2 \, 3)$. The group $A$ is isomorphic to the dihedral group $D_4$ of order $8$. Define $b \in \St(1)$ by
\[
	b = (0: b,\, 1: s,\, 2: s, \, 3: sr),
\]
and $B = \la b \ra$. Let $G$ be the \csp group defined by $A$ and $B$.

Clearly, $G$ is not orbitwise-abelian, since the rotation $r$ acts transitively on $X$. Thus $G$ is in particular not strongly orbitwise-abelian and does not satisfy the hypothesis of \cref{Theorem}{thm:A}.

To prove that $G$ does not satisfy the hypothesis of \cref{Theorem}{thm:B}, we have to consider all generating sets of $B$. But $B$ is cyclic of order two, whence it suffices to determine $\Sigma_{\{b\}}$. The point stabiliser of $0$ in $A$ is generated by $s$, hence the set $\mathfrak C(a, x)$ is of cardinality at most two for all $a \in A, x \in X$. For our argument, it suffices to calculate
\[
	\mathfrak C(r^{\pm 1}, x) = \{ r^{\pm 1} \} \quad \text{ for all } x \in X
	\quad\text{and}\quad
	\mathfrak C(sr, 0) = \{ sr, sr^3 \}.
\]
This yields
\[
	\sigma_{\{b\}}(r^{\pm 1}, x) = \la sr \ra \cdot \la s, sr \ra' = \la sr \ra \times \la r^2 \ra,
\]
for each $x \in X$, hence $\Sigma_{\{b\}}(r^{\pm 1}) = \la sr \ra \times \la r^2 \ra$, and
\[
	\sigma_{\{b\}}(sr, 0) = \la s, sr \ra \cdot \la s, sr \ra' = A,
\]
whence $sr \in \Sigma_{\{b\}}(sr)$. Thus $\Sigma_{\{b\}}$ is not eventually trival.

We now prove that $G$ is periodic. To achive this, we use \cref{Proposition}{prop:torsion_by_length} with the set
\[
	T = B \cup{}^GA \cup {}^B\left(\la bsr \ra \cdot \{ 1_A, bs \}\right).
\]
A standard computation shows that the elements of $T$ have finite order: The first two sets in the union clearly consist of elements of finite order. For the third subset, first calculate
\begin{align*}
	(bs)^2 &= (1_A, r, 1_A, r^3), &&\text{ hence } \ord(bs) = 8, \text{ and}\\
	(bsr)^{2} &= (bs, (bs)^{-1}, r, r^3), &&\text{ hence } \ord(bsr) = 16.
\end{align*}
The first equation also implies that $(bs)^{2n}r^{\pm 1}$ is of finite order for all $n \in \N$. Similarly
\begin{align*}
	((bs)^{2n+1}r^{2})^2 &= ((b, sr^{-n}, s, sr^{n+1})sr^2)^2\\
	&= (bs, 1_A, (bs)^{-1}, 1_A)
\end{align*}
has finite order for all $n \in \N$. Using this, we see that for arbitary $n \in \N$
\begin{align*}
	((bsr)^{2n} bs)^2 &= \bigg(\left((bs)^n, (bs)^{-n}, r^n, r^{-n}\right) bs\bigg)^2\\
	&= \bigg(\left((bs)^nb, (bs)^{-n}s, sr^{-n}, sr^{n+1}\right) s\bigg)^2\\
	&= (1_A, (bs)^{-n}r^{n+1}, 1_A, ((bs)^{-n}r^{n+1})^{-1})
\end{align*}
is of finite order. Finally,
\begin{align*}
	((bsr)^{2n + 1} bs)^4 &= ((bs, (bs)^{-1}, r, r^3)^n b(^{sr}b)r^3)^4\\
	&= (((bs)^{n+1}, (bs)^{-n-1}, r^{n+1}, r^{-n-1})r^3)^4 = 1_A.
\end{align*}
Thus the elements of $\la bsr \ra \cdot \{ 1_A, bs \}$ are of finite order, and $T$ is a valid choice for (\ref{eq:lr}).

Now let $g = {}^{a_0}b\;^{a_1}b \dots {}^{a_{n-1}}b\; a_n \in G$ be an element of syllable length $n \in \N$. If $x \in X$ is a letter such that $\syl(g\|_x) = n$, we have $g\|_x|^\epsilon \in \Sigma_{\{b\}}(a_n)$, as we have seen in the proof of \cref{Theorem}{thm:B}. Thus if $a_n$ is of order four, i.e.\ $a_n \in \{r, r^3\}$, our computation above shows that $g\|_x|^\epsilon \in \Sigma_{\{b\}}(r^{\pm 1}) = \la sr \ra \times \la r^2 \ra$ is of order two. Thus it is enough to prove that for every $g$ such that $g|^\epsilon$ has order two, there is a number $k \in \N$ such that $g\|_u \in T$ or $\syl(g\|_u) < n$ for all $u \in X^k$.

Assume that $a_n$ is of order two and that $\syl(g\|_x) = n$. Then we calculate $g\|_{x}$ as we have done in \cref{Lemma}{lem:non-reduction in one step}. If $n$ is even, we see that either
\begin{align*}
	\begin{array}{rcccccccccc}
		g\|_{x} = && b & \cdot & c_0 & \cdot & b & \cdots & c_{\frac n 2} & \cdot & b\\
		& \cdot & c_{\frac n 2 + 1} & \cdot & b & \cdot & c_{\frac n 2 + 2} & \cdots & b & \cdot & c_{n-1},
	\end{array}
\end{align*}
where $c_i \in \{b|_x \mid x \in \dot X\} = \{ s, sr \}$ for $i \in \{ 0, \dots, n-1 \}$, or $g\|_x$ is conjugate to such an expression. If $n$ is odd, we obtain
\begin{align*}
	\begin{array}{rcccccccccc}
		g\|_{x} = && b & \cdot & c_0 & \cdot & b & \cdots & b & \cdot & c_{\frac {n+1}2}\\
		& \cdot & c_{\frac {n+1}2 + 1} & \cdot & b & \cdot & c_{\frac{n+1}2 + 2} & \cdots & c_{n-1} & \cdot & b,
	\end{array}
\end{align*}
with $c_i \in \{ s, sr \}$ for $i \in \{ 0, \dots, n-1 \}$.
In this case, the decomposition beginning with a $c_i$-letter is of syllable length $n-2$, since the two instances of $b$ in the middle of the word cancel each other. Assume that there is a letter $c_i$, for any $i \neq n - 1$ in case $n$ is even, and $i \notin \{(n+1)/2, (n+1)/2+1\}$ if $n$ is odd, such that $c_i = s$. Then there are elements $\hat a_i \in A$ for $i \in \{0, \dots, n\}$ such that we may represent
\[
	g\|_{x} = {}^{\hat a_0}b\,^{\hat a_1}b \dots {}^{\hat a_{n-1}}b\,\hat a_n
\]
and $\hat a_{i+1} = \hat a_i s$. Assume that $\syl(g\|_{x\star y}) = n$ for some letter $y \in X$. According to \cref{Lemma}{lem:non-reduction in one step}, we calculate
\[
	g\|_{x\star y} = \prod_{j = 0}^{\ell_{g\|_{x}}(y)} \prod_{i = 0}^{n-1} b|_{x_1.\hat a_n^j \hat a_i}.
\]
Since $\syl(g\|_{x\star y}) = n$ there is some $j \in \{0, \dots, \ell_{g\|_{x}}(y)\}$ such that $y.\hat a_n^j a_i = 0$. But then $y.\hat a_n^j a_{i+1} = y.\hat a_n^j a_{i}s = 0.s = 0$, and two of the sections that evaluate to $b$ cancel each other. Thus, $\syl(g\|_{x\star y}) \leq n - 2$, which is a contradiction. Thus either the syllable length of $g$ reduces when taking stabilised sections at words in $X^2$, or
\[
	g\|_x \in \{(bsr)^{n-1} b \cdot \{ s, sr \} \mid n \text{ even} \} \cup \{(bsr)^{\frac{n-1}2}b \cdot \{ 1_A, r^3 \} \cdot (bsr)^{\frac{n-1}2}b \mid n \text{ odd} \}.
\]
But this set is contained in $T$, hence the conditions of \cref{Proposition}{prop:torsion_by_length} are satisfied and $G$ is periodic.



\end{document}